\documentclass{amsart}
\usepackage{amsthm}
\usepackage{amsmath}
\usepackage{amsfonts}
\usepackage{amssymb}
\usepackage{graphicx}
\usepackage{float}
\usepackage{rotating}
\restylefloat{figure}
\usepackage{mathrsfs}
\usepackage{fancyhdr}
\usepackage[breaklinks]{hyperref}
\usepackage[numbers,sort&compress]{natbib}
\theoremstyle{plain}
\newtheorem{lem}{Lemma}

\newtheorem{thm}[lem]{Theorem}

\theoremstyle{definition}

\newcommand{\R}{\mathbb R}

\newcommand{\N}{\mathbb N}
\newcommand{\dx}{\,dx}
\newcommand{\dz}{\,dz}
\renewcommand{\d}{\,d}

\newcommand{\eps}{\varepsilon}
\expandafter\mathchardef\expandafter\varphi\number\expandafter\phi\expandafter\relax
\expandafter\mathchardef\expandafter\phi\number\varphi
\newcommand{\norm}[1]{\left\Vert#1\right\Vert}
\newcommand{\ska}[2]{\left\langle #1,#2\right\rangle}
\renewcommand{\autoref}[1]{\text{Eq.}~\eqref{#1}}
\newcommand{\bea}{\begin{eqnarray}}
\newcommand{\eea}{\end{eqnarray}}
\newcommand{\beq}{\begin{equation}}
\newcommand{\eeq}{\end{equation}}
\begin{document}
\title{The abstract quasilinear Cauchy problem for a MEMS model with two free boundaries}
\author{Martin Kohlmann}
\address{Dr.\ Martin Kohlmann, Goerdelerstraße 36, 38228 Salzgitter, Germany}
\email{martin\_kohlmann@web.de}
\keywords{MEMS, free boundary problem, local and global well-posedness, non-existence, asymptotic stability, small aspect ratio limit}
\subjclass[2010]{35R35, 35M33, 35B30, 35Q74, 74M05}
\date{\today}
\begin{abstract} In this paper, we reformulate a mathematical model for the dynamics of an idealized electrostatically actuated MEMS device with two elastic membranes as an initial value problem for an abstract quasilinear evolution equation. Applying the Contraction Mapping Theorem, it is shown that the model is locally well-posed in time for any value of the source voltage of the device. In addition it is proven that the MEMS model considered here possesses global solutions for small source voltages whereas for large source voltages solutions of the model have a finite maximal existence time. Furthermore, we comment on the relationship of our model to its stationary version and to its small aspect ratio limit by showing that there exists a unique exponentially stable steady state and by proving convergence towards a solution of the narrow gap model in the vanishing aspect ratio limit.\ Our results extend the discussion of the elliptic-parabolic MEMS model presented in \cite{mk13_2} leading to a Cauchy problem for a semilinear abstract evolution equation.
\end{abstract}
\maketitle
\section{Introduction and main results}\label{sec_intro}
Micro-electro mechanical systems (MEMS) are small devices that operate on the principle of electrostatic actuation: applying a potential difference between certain mechanical components of the device causes an electric field and hence a Coulomb force resulting in a mechanical deformation. There are are wide range of applications of MEMS to report on: MEMS are used as microsensors and microactuators, they appear as components of accelerometers and gyroscopes, they have commercial applications, e.g., in microphones and mobile phones, and Bio-MEMS are used in medical and health technology. In recent years, MEMS have also become a flourishing field of research in mathematics as various types of models for such devices have been proposed. Most often these models are concerned with an idealized device consisting of a deflecting membrane suspended above a fixed ground plate.

\begin{figure}[H]
\begin{center}
\includegraphics[width=9cm]{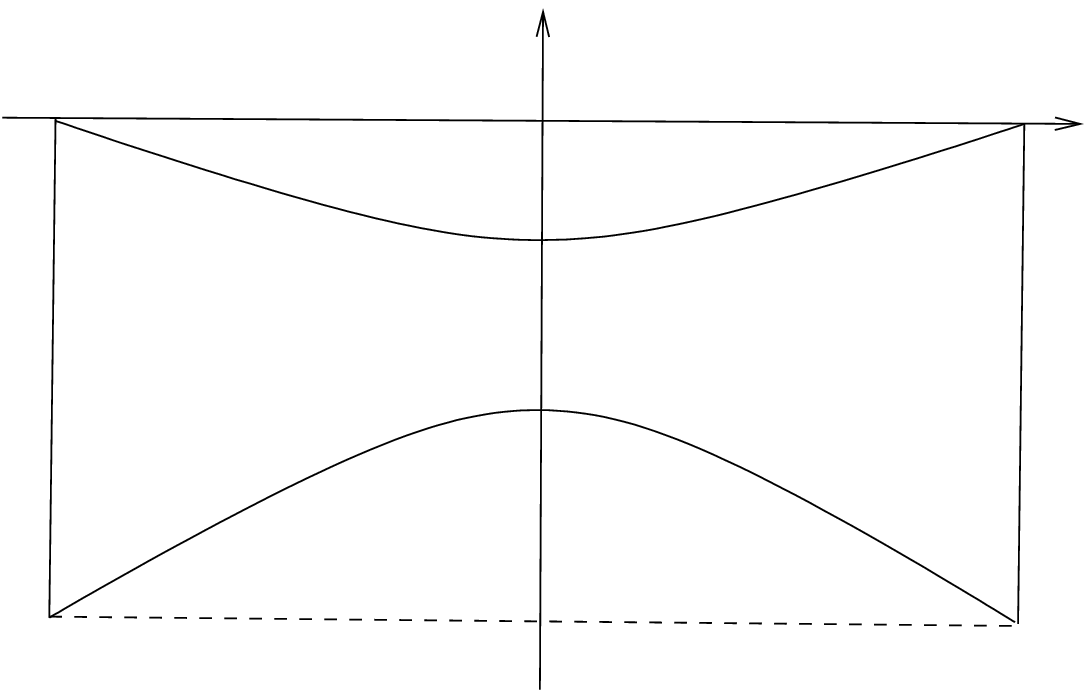}
\end{center}
\caption{An idealized model for an electrostatic MEMS device with two free boundaries.}
\label{figabstr}
\begin{picture}(0,0)
\put(-60,95){$\Omega_{u,v}$}
\put(-90,55){$\Gamma_v$}
\put(-90,132){$\Gamma_u$}
\put(-46,150){\vector(0,-1){22}}
\put(-46,33){\vector(0,1){37}}
\put(-40,135){$u(x)$}
\put(-40,45){$v(x)+1$}
\put(10,153){$z=0$}
\put(10,35){$z=-1$}
\put(108,153){$+1$}
\put(-122,153){$-1$}
\put(3,170){$z$}
\put(130,147){$x$}
\end{picture}
\end{figure}

In the paper at hand we discuss a moving boundary problem for a so-called DFM device, i.e., a MEMS with double freestanding membranes as explained in, e.g., \cite{FTYS13}. Our model can be derived as follows: Let $H,L>0$ and denote by $(\hat x,\hat z)$ coordinates of the two-dimensional rectangular domain $R=(-L,L)\times(0,-H)$. We consider two thin, conductive and elastic membranes of length $2L$ and distance $H$ located at the upper and the lower boundary of $R$ which should be held fixed at $(\pm L,0)$ and $(\pm L,-H)$. Moreover, we assume that the permittivity of the medium filling the interior of $R$ is equal to one. When a non-zero source voltage $V$ is applied to the device, an electric field sets up causing a deformation of the membranes whose displacements are then modeled by functions $\hat u$, $\hat v$ satisfying $-H<\hat v(\hat x)<\hat u(\hat x)<0$, for $\hat x\in (-L,L)$, and $(\hat u,\hat v)(\pm L)=(0,-H)$. Let $\hat\phi(\hat x,\hat z)$ denote the electrostatic potential defined in the region
$$\hat\Omega_{\hat u,\hat v}:=\{(\hat x,\hat z)\in R;\,\hat v(\hat x)<\hat z<\hat u(\hat x)\}$$
between the membranes. Then $\hat\phi$ is a solution to the Laplace equation,
$$\partial_{\hat x}^2\hat\phi+\partial_{\hat z}^2\hat\phi=0\quad\text{in }\hat\Omega_{\hat u,\hat v},$$
with the boundary conditions
\begin{align}
\hat\phi(\hat x,\hat v(\hat x))&=0,&\hspace{-3cm}\hat x\in(-L,L),\nonumber\\
\hat\phi(\hat x,\hat u(\hat x))&=V,&\hspace{-3cm}\hat x\in(-L,L).\nonumber
\end{align}
We assume that the continuous extension of $\hat\phi$ to the lateral boundary of $R$ depends linearly on $\hat z$. The total potential energy $E(\hat u,\hat v)$ of the device is the sum of the electrostatic energy determined by the square of the gradient of the potential plus the elastic energy determined by the change of the length of the elastic membranes. To be able to compare the strengths of the mechanical and electrical forces in the device, we also introduce surface tension coefficients $T_1,T_2>0$ so that
\begin{align} E(\hat u,\hat v) & =\frac{\eps_0}{2}\int_{-L}^L\int_{\hat v(\hat x)}^{\hat u(\hat x)}|\nabla\hat\phi(\hat x,\hat z)|^2\d\hat x d\hat z
+T_1\int_{-L}^L\left(\sqrt{1+(\partial_{\hat x}\hat u(\hat x))^2}-1\right)d\hat x \nonumber\\
&\qquad + T_2\int_{-L}^L\left(\sqrt{1+(\partial_{\hat x}\hat v(\hat x))^2}-1\right)d\hat x, \nonumber
\end{align}
where $\eps_0$ is the permittivity of free space. We now define dimensionless variables
$$x=\frac{\hat x}{L},\quad z=\frac{\hat z}{H},\quad u=\frac{\hat u}{H},\quad v=\frac{\hat v}{H},\quad\phi=\frac{\hat\phi}{V}$$
and parameters
$$\eps=\frac{H}{L},\quad\lambda = \frac{\eps_0V^2}{2\eps^3T_1L},\quad\mu=\frac{\eps_0V^2}{2\eps^3T_2L},$$
we introduce the sets
\begin{align}
\Omega_{u,v}&=\{(x,z)\in(-1,1)\times(-1,0);\,v(x)<z<u(x)\}, \nonumber\\
\Gamma_u&=\{(x,u(x));\,x\in I\}, \nonumber\\
\Gamma_v&=\{(x,v(x));\,x\in I\}, \nonumber
\end{align}
shown in Figure~\ref{figabstr}, and $I=(-1,1)$ and define the operators $\nabla_\eps=(\eps\partial_x,\partial_z)$ and $\Delta_\eps=\eps^2\partial_x^2+\partial_z^2$.
The total energy of the device can be rewritten as
\begin{align}
E(\hat u,\hat v)&=\frac{\eps_0}{2\eps}V^2\int_{\Omega_{u,v}}|\nabla_\eps\phi(x,z)|^2\dx dz+T_1L\int_I\left(\sqrt{1+\eps^2(\partial_xu)^2}-1\right)dx \nonumber\\
&\qquad +T_2L\int_I\left(\sqrt{1+\eps^2(\partial_xv)^2}-1\right)dx \nonumber\\
& = \int_{I}\mathcal L\left(x,u,v,\partial_xu,\partial_xv\right)dx\nonumber
\end{align}
with $\mathcal L$ denoting the Lagrangian density. The Euler-Lagrange equations $\partial_x\partial_{\partial_xu}\mathcal L-\partial_u\mathcal L=0$ and $\partial_x\partial_{\partial_xv}\mathcal L-\partial_v\mathcal L=0$ take the form
\begin{align}
\partial_x\left(\frac{\partial_xu}{\sqrt{1+\eps^2(\partial_xu)^2}}\right)-\lambda|\nabla_\eps\phi(x,u(x))|^2&=0, \nonumber\\
\partial_x\left(\frac{\partial_xv}{\sqrt{1+\eps^2(\partial_xv)^2}}\right)+\mu|\nabla_\eps\phi(x,v(x))|^2&=0. \nonumber
\end{align}
We now assume that $u$ and $v$ also depend on time $\hat t$. Then $\partial_{\hat t}^2u$ models the acceleration of $\Gamma_u$ and $\Gamma_v$ in the associated evolution problem. Regarding the left-hand sides of the Euler-Lagrange equations above as forces on $\Gamma_u$ and $\Gamma_v$ and considering a damping force that is proportional to the velocity $\partial_{\hat t}u$, Newton's Second Law yields that
\begin{align}
\rho_1\delta_1\partial_{\hat t}^2u+a\partial_{\hat t}u&=\partial_x\left(\frac{\partial_xu}{\sqrt{1+\eps^2(\partial_xu)^2}}\right)-\lambda|\nabla_\eps\phi(x,u(x))|^2, \nonumber\\
\rho_2\delta_2\partial_{\hat t}^2v+a\partial_{\hat t}v&=\partial_x\left(\frac{\partial_xv}{\sqrt{1+\eps^2(\partial_xv)^2}}\right)+\mu|\nabla_\eps\phi(x,v(x))|^2, \nonumber
\end{align}
where $\rho_1,\rho_2$ and $\delta_1,\delta_2$ denote the mass density per unit volume of the membranes and the membrane thicknesses respectively and $a$ is a damping constant. With
$$t=\frac{\hat t}{a},\quad\gamma_1=\frac{\sqrt{\rho_1\delta_1}}{a},\quad\gamma_2=\frac{\sqrt{\rho_2\delta_2}}{a}$$
we arrive at the non-dimensionalized equations
\begin{align}
\gamma_1^2\partial_t^2u+\partial_tu & = \partial_x\left(\frac{\partial_xu}{\sqrt{1+\eps^2(\partial_xu)^2}}\right)-\lambda|\nabla_\eps\phi(x,u(x))|^2, \nonumber\\
\gamma_2^2\partial_t^2v+\partial_tv & = \partial_x\left(\frac{\partial_xv}{\sqrt{1+\eps^2(\partial_xv)^2}}\right)+\mu|\nabla_\eps\phi(x,v(x))|^2. \nonumber
\end{align}
In this paper, we will assume that $\gamma_1,\gamma_2\ll 1$ meaning that the damping forces dominate over the inertial forces. Given initial values $u_0$ and $v_0$ for the functions $u$ and $v$, we thus discuss the following system of equations:
\begin{align}
-\Delta_\eps\phi & = 0, & \text{in }\Omega_{u,v},\,t>0,\label{originalproblem1}\\
\phi & = \frac{z-v}{u-v}, & \text{on }\partial\Omega_{u,v},\,t>0,\label{originalproblem2}\\
\partial_tu-\partial_x\left(\frac{\partial_xu}{\sqrt{1+\eps^2(\partial_xu)^2}}\right) & = -\lambda|\nabla_\eps\phi|^2,&\text{on }\Gamma_u, \,t>0,\label{originalproblem3}\\
\partial_tv-\partial_x\left(\frac{\partial_xv}{\sqrt{1+\eps^2(\partial_xv)^2}}\right) & = \mu|\nabla_\eps\phi|^2,&\text{on }\Gamma_v,\,t>0,\label{originalproblem4}\\
u(t,\pm 1) & = 0,&t>0,\label{originalproblem5}\\
v(t,\pm 1) & =-1,&t>0,\label{originalproblem6}\\
u(0,x) & = u_0,&x\in I,\label{originalproblem7}\\
v(0,x) & = v_0,&x\in I.\label{originalproblem8}
\end{align}
Note that \eqref{originalproblem1}--\eqref{originalproblem8} is a free boundary problem as the domain $\Omega_{u,v}$ and its boundary components $\Gamma_u,\Gamma_v$ have to be determined together with the solution $(u,v,\phi)$.
Several simplified models of \eqref{originalproblem1}--\eqref{originalproblem8} have been studied recently: In \cite{mk13_2} we have assumed that the deformation of the membranes is small so that, in the equations on the free boundaries, the curvature terms on the left-hand sides of \eqref{originalproblem3}--\eqref{originalproblem4} can be replaced by the linear terms $-\partial_x^2u$ and $-\partial_x^2v$. In this case, the evolution of the membranes is described by two heat equations with a right-hand side proportional to the square of the gradient of the potential on the boundary. In \cite{mk13} the stationary version of the MEMS model with two free boundaries and linear stretching terms has been discussed. For $v\equiv -1$, the problem \eqref{originalproblem1}--\eqref{originalproblem8} models the evolution of a free membrane suspended above a fixed ground plate. Various analytical results on this type of a MEMS have been obtained in recent years: \cite{ELW12,FMPS0607,GG08,G08a,G08b,H11,LW13I,LW13II} refer to the parabolic problem, \cite{CFT14,FMP03,G10,KLNT11,LW13I} discuss the problem with a hyperbolic evolution equation and in \cite{C07,P0102,PT00,P01,WL12} the stationary model is presented. The corresponding model with an additional curvature term is discussed in \cite{ELW13} and our derivation of \eqref{originalproblem1}--\eqref{originalproblem8} refines Lauren\c{c}ot's line of arguments therein.

For $\eps\to 0$, one obtains the so-called small aspect ratio model from \eqref{originalproblem1}--\eqref{originalproblem8}:
\begin{align}
\phi & = \frac{z-v}{u-v}, & \text{in }{\Omega}_{u,v}\cup\partial\Omega_{u,v},\,t>0,\label{sarm1}\\
\partial_tu-\partial_x^2u & = -\frac{\lambda}{(u-v)^2},&x\in I,\, t>0,\label{sarm2}\\
\partial_tv-\partial_x^2v & = \frac{\mu}{(u-v)^2},&x\in I,\, t>0,\label{sarm3}\\
u(t,\pm 1) & = 0,&t>0,\label{sarm4}\\
v(t,\pm 1) & =-1,&t>0,\label{sarm5}\\
u(0,x) & = u_0,&x\in I,\label{sarm6}\\
v(0,x) & = v_0,&x\in I.\label{sarm7}
\end{align}
The problem \eqref{sarm1}--\eqref{sarm7} already appeared in \cite{mk13_2} where we proved that solutions of the MEMS model of \cite{mk13_2} with $\eps>0$ converge towards solutions of \eqref{sarm1}--\eqref{sarm7} in the vanishing aspect ratio limit. The small aspect ratio limit of the MEMS model with a fixed ground plate is a subject of \cite{BP11,ELW13,FMP03,FMPS0607,EGG10,GG0607,GG08,H11,GHW09,G08a,G08b,LY07,PB03}.

The plan of the present paper is to apply and refine the chain of arguments used in \cite{ELW13,mk13_2} in order to obtain results on solutions of \eqref{originalproblem1}--\eqref{originalproblem8} where we have to cope with additional curvature terms compared to the model in \cite{mk13_2}. In doing so, our first aim is to show that \eqref{originalproblem1}--\eqref{originalproblem8} possesses a unique maximal solution for any pair of values $(\lambda,\mu)$. To this end, we solve the elliptic problem \eqref{originalproblem1}--\eqref{originalproblem2} for the potential and then rewrite the system \eqref{originalproblem3}--\eqref{originalproblem8} as an initial value problem for an abstract quasilinear evolution equation whose solution is obtained from the variation of constants formula and the Contraction Mapping Theorem. Our main endeavour is to prove the Lipschitz continuity of the right-hand side with respect to the topology of $W^{2-\xi}_q(I)\times W^{2-\xi}_q(I)$, $\xi>0$; see also the semilinear problems in \cite{mk13_2,ELW12} where this has been achieved for $\xi=0$. Our first main result which is the analog of \cite[Theorem 2]{mk13_2} and \cite[Theorem~1.1]{ELW13} reads as follows.
\thm\label{thm_lwp} Let $q\in(2,\infty)$ and $\eps\in(0,1)$ and consider initial values $(u_0,v_0)\in W_{q}^2(I)\times W_{q}^2(I)$ satisfying $(u_0,v_0)(\pm 1)=(0,-1)$ and $-1\leq v_0(x)<u_0(x)\leq 0$ for all $x\in I$. Then:
\begin{itemize}
\item[(i)] There exists $r>0$ such that for $\norm{u_0}_{W^2_{q}(I)}, \norm{v_0+1}_{W^2_{q}(I)}<r$ and for any $\lambda,\mu>0$, there is a unique maximal solution $(u,v,\phi)(t)$, $t\in[0,T_\eps)$, $T_\eps>0$, to \eqref{originalproblem1}--\eqref{originalproblem8} with regularity
$$u,v\in C([0,T_\eps),W^2_{q}(I))\cap C^1([0,T_\eps),L_q(I)),\,\phi\in W_2^2(\Omega_{u(t),v(t)})$$
so that $-1\leq v<u\leq 0$ on $[0,T_\eps)\times I$.
\item[(ii)] If for each $\tau>0$ there exists $\kappa(\tau)\in(0,\tfrac{1}{2})$ such that $u(t)-v(t)\geq 2\kappa(\tau)$ and $\norm{u(t)}_{W_q^2(I)},\norm{v(t)+1}_{W_q^2(I)}\leq\kappa(\tau)^{-1}$ for $t\in[0,T_\eps)\cap[0,\tau]$, then the solution exists globally in time, i.e., $T_\eps=\infty$.
\item[(iii)] If $u_0$ and $v_0$ are even functions on $I$, then $(u,v,\phi)$ is even in $x$ on $[0,T_\eps)\times I$.
\item[(iv)] Given $\kappa\in(0,\tfrac{1}{2})$, there exist $m_0(\kappa),r_0(\kappa)>0$ such that, for $\max\{\lambda,\mu\}<m_0(\kappa)$ and $\norm{u_0}_{W^2_q(I)},\norm{v_0+1}_{W^2_q(I)}<r_0(\kappa)$, one has $T_\eps=\infty$, $u(t)-v(t)\geq 2\kappa$ and $\norm{u(t)}_{W_q^2(I)}$ and $\norm{v(t)+1}_{W_q^2(I)}$ are bounded by a positive constant only depending on $\kappa$.
\end{itemize}
\endthm\rm
Observe that, in contrast to Theorem~2 of \cite{mk13_2} and Theorem~1.1 of \cite{ELW13}, we have to assume that the initial values $(u_0,v_0)$ are sufficiently small in $W^2_q(I)\times W^2_q(I)$ here. A proof of Theorem~\ref{thm_lwp} can be found in Section~\ref{sec_lwp}. The methods used in Section~\ref{sec_lwp} also yield that solutions of \eqref{originalproblem1}--\eqref{originalproblem8} converge towards a solution of \eqref{sarm1}--\eqref{sarm7} for $\eps\to 0$. We present a proof of the following theorem which is the analog of \cite[Theorem 10]{mk13_2} and \cite[Theorem~1.4]{ELW13} and which justifies rigorously the relationship between the original problem and its small aspect ratio limit. Here, $\mathbf{1}_{A}$ denotes the indicator function of the set $A\subset\R^2$.
\thm\label{thm_sar} Let $\lambda,\mu>0$, $q\in(2,\infty)$ and let $(u_0,v_0)\in W^2_q(I)\times W^2_q(I)$ satisfying the assumptions in Theorem~\ref{thm_lwp} be given. For $\eps\in(0,1)$, the unique solution to \eqref{originalproblem1}--\eqref{originalproblem8} with initial values $(u_0,v_0)$ and the maximal interval of existence $[0,T_\eps)$ is denoted by $(u_\eps,v_\eps,\phi_\eps)(t)$. Then there are $\tau>0$, $\eps_*\in(0,1)$ and $\kappa_1\in(0,\tfrac{1}{2})$ 
such that $T_\eps\geq\tau$, $u_\eps(t)-v_\eps(t)\geq 2\kappa_1$ and $\norm{u_\eps(t)}_{W_q^2(I)},\norm{v_\eps(t)+1}_{W_q^2(I)}\leq\kappa_1^{-1}$, for all $(t,\eps)\in[0,\tau]\times(0,\eps_*)$. Moreover, the small aspect ratio model \eqref{sarm1}--\eqref{sarm7} has a unique solution $(u_*,v_*,\phi_*)$ satisfying
$$u_*,v_*\in C([0,\tau],W^2_{q}(I))\cap C^1([0,\tau],L_q(I)),$$
$-1\leq v_*(t)<u_*(t)\leq 0$, $u_*(t)-v_*(t)\geq 2\kappa_1$, for all $t\in[0,\tau]$, and there is a null sequence $(\eps_n)_{n\in\N}\subset(0,\eps_*)$ such that
\begin{align}
(u_{\eps_n},v_{\eps_n})&\to (u_*,v_*) &\text{in }C^{1-\theta}([0,\tau],W_q^{2\theta}(I)),\,\theta\in(0,1), \nonumber\\
\phi_{\eps_n}(t)\mathbf{1}_{\Omega_{u_{\eps_n}(t),v_{\eps_n}(t)}}&\to\phi_*(t)\mathbf{1}_{\Omega_{u_*(t),v_*(t)}} &\text{ in }L_2(I\times(-1,0)),\,t\in[0,\tau], \nonumber
\end{align}
as $n\to\infty$. Furthermore, there is $\Lambda(\kappa)>0$ such that, for $\lambda,\mu<\Lambda(\kappa)$, the statements of the theorem hold true for any $\tau>0$.
\endthm\rm
In particular, Theorem~\ref{thm_sar} guarantees that the maximal existence times $T_\eps$ are bounded from below when sending $\eps\to 0$. Again, in contrast to the models discussed in \cite{mk13_2} and \cite{ELW13}, an additional condition on the norm of the initial values occurs in the above theorem.

The effectiveness of our MEMS device is limited when increasing the source voltage as the membranes might come close and closer and finally touch. This phenomenon is called pull-in stability and has already been discussed for related models, see, e.g., \cite{BGP00,ELW12,ELW13,ELW13b,FMP03,FMPS0607,G08b,P01,WL12}. It is plausible to expect that for small voltage values the problem \eqref{originalproblem1}--\eqref{originalproblem8} has a global solution and that for $\lambda$ and $\mu$ sufficiently large, there is no steady state of \eqref{originalproblem1}--\eqref{originalproblem8}. Recall that Theorem~\ref{thm_lwp}.(iv) implies that solutions $(u,v,\phi)(t)$ exist globally in time in the sense that neither touchdown of the membranes nor blow up of the displacements in $W^2_q(I)\times W^2_q(I)$ occurs, provided $\lambda$ and $\mu$ and the initial values are sufficiently small. Next, we complement this result by a non-existence theorem for high voltages.\ We will concentrate on displacements $u$ and $v$ that have a positive distance to $\{z=-1\}$ and $\{z=0\}$ respectively, as touchdown on $\{z=-1\}$ or $\{z=0\}$ is reminiscent of the associated MEMS problem with only one free membrane. For sufficiently large values $\lambda$ and $\mu$, we divine that $T_\eps<\infty$ and
\beq
\limsup_{t\to T_\eps}\norm{(u,v)(t)}_{W^2_q(I)\times W^2_q(I)}=\infty\quad\text{or}\quad\liminf_{t\to T_\eps}\min\{u(t)-v(t)\}=0.\label{condTeps}
\eeq
It will remain an open problem whether the membranes certainly smash-up when $T_\eps<\infty$. By \eqref{condTeps}, the displacements might also blow up in $W^2_q(I)\times W^2_q(I)$ contradicting the physical expectation that there is collision of the membranes in the interior of the device for finite maximal existence times. A similar ambiguity has been observed in \cite{LW13I}.
In Section~\ref{sec_gwp} we present proofs of the following theorems.
\thm\label{thm_nonex1} Let $q\in(2,\infty)$ and $\eps\in(0,1)$. There exists a positive number $\xi_0(\eps)$ such that for $\max\{\lambda,\mu\}>\xi_0(\eps)$ the stationary problem \eqref{originalproblem1}--\eqref{originalproblem8} possesses no steady state solution $(u,v,\phi)$ with $u,v\in W^2_q(I)$ and $\phi\in W_2^2(\Omega_{u,v})$ satisfying $-1\leq v<u\leq 0$ on $I$. In addition $\xi_0(\eps)\to 2$ for $\eps\to 0$.
\endthm\rm
\thm\label{thm_nonex2} Let $q\in(2,\infty)$ and $u_0,v_0\in W^2_q(I)$ satisfying the assumptions in Theorem~\ref{thm_lwp} and the additional assumption $(u_0,v_0)(-x)=(u_0,v_0)(x)$, for all $x\in I$, be given. There exists $\eps_0\in(0,1)$ such that for all $\eps\in(0,\eps_0)$ the following holds true: If $\max\{\lambda,\mu\}>4/\eps$, the displacements $(u,v)$ do not blow up in $W^2_q(I)\times W^2_q(I)$ and $v,-u-1\leq c$, for some $c<0$, then the maximal existence time of the solution $(u,v,\phi)$ obtained in Theorem~\ref{thm_lwp} is finite and
\beq\label{Tepsabsch}T_\eps\leq\frac{\norm{u_0-v_0}_{L_1(I)}}{\max\{\lambda,\mu\}-4/\eps}.\eeq
If equality holds in \eqref{Tepsabsch}, there is touchdown of the membranes in the sense that $\liminf_{t\to T_\eps}\min\{u(t)-v(t)\}=0$.
\endthm\rm
Note that Theorem~\ref{thm_nonex1} and Theorem~\ref{thm_nonex2} are the analogs of \cite[Theorem 5]{mk13_2}, \cite[Theorem~1.3]{ELW13} and \cite[Theorem 3]{ELW13b}.

Finally, in Section~\ref{sec_steady}, it will be established that, for any $\kappa\in(0,1/2)$, \eqref{originalproblem1}--\eqref{originalproblem8} possesses a unique steady state so that the boundary components have distance at least $2\kappa$ and the $W_q^2(I)$-norms of the first and second component are bounded by $\kappa^{-1}$. Moreover exponential stability of this steady state is shown using the Principle of Linearized Stability. The following theorem is the analog of \cite[Theorem 6]{mk13_2} and \cite[Theorem~1.2]{ELW13}.
\thm\label{thm_stability} Let $q\in(2,\infty)$, $\eps\in(0,1)$ and $\kappa\in(0,\tfrac{1}{2})$ be fixed.
\begin{enumerate}
\item[(i)] There are $\delta(\kappa)>0$ and and analytic function $[0,\delta)^2\to W^2_{q}(I)\times W_{q}^2(I)$, $\Lambda\to U_\Lambda=(U_{\Lambda,1},U_{\Lambda,2})$, such that $(U_\Lambda,\Phi_\Lambda)$ is for each $\Lambda=(\lambda,\mu)\in(0,\delta)^2$ the unique steady state of \eqref{originalproblem1}--\eqref{originalproblem8} satisfying $U_{\Lambda,1}-U_{\Lambda,2}\geq 2\kappa$ and $\norm{U_{\Lambda,1}}_{W^2_{q}(I)},\norm{U_{\Lambda,2}+1}_{W^2_{q}(I)}\leq\kappa^{-1}$ and $\Phi_\Lambda\in W_2^2(\Omega_{U_{\Lambda,1},U_{\Lambda,2}})$ is the potential associated with $U_\Lambda$. Moreover, $U_{\Lambda,1}$ and $-U_{\Lambda,2}$ are convex and even with $U_{(0,0)}=(0,0)$ and $x\mapsto\Phi_\Lambda(x,z)$ is even on $I$.
\item[(ii)] Let $\Lambda\in(0,\delta)^2$. There are numbers $\omega_0,\varrho,R>0$ such that for each pair of initial values $u_0,v_0\in W_{q}^2(I)$ satisfying the assumptions in Theorem~\ref{thm_lwp} and the additional assumption $\norm{(u_0,v_0)-U_{\Lambda}}_{W_{q}^2(I)\times W_{q}^2(I)}<\varrho,$
    the associated solution $(u,v,\phi)$ to \eqref{originalproblem1}--\eqref{originalproblem8} exists globally in time with $u(t)-v(t)>0$ and
    \begin{align}
    &\norm{(u,v)-U_\Lambda}_{W_{q}^2(I)\times W_{q}^2(I)}+\norm{(u_t,v_t)}_{L_q(I)\times L_q(I)}\nonumber\\
    &\hspace{2cm}\leq Re^{-\omega_0t}\norm{(u_0,v_0)-U_{\Lambda}}_{W_{q}^2(I)\times W_{q}^2(I)},\quad\forall t\geq 0.\nonumber
    \end{align}
\end{enumerate}
\endthm\rm
A convergence result similar to Theorem~\ref{thm_stability}.(ii) holds true for the first component $\phi$ of the solution, cf. Section~\ref{sec_steady} for the technical details.
\section{Local and global well-posedness and the small aspect ratio limit}\label{sec_lwp}
In this section, we present proofs of Theorem~\ref{thm_lwp} and Theorem~\ref{thm_sar}.
Let us first introduce our notation and recall some preliminary results: Let $\Omega:=I\times(0,1)$ and consider the time-dependent transformation of coordinates $T=T_{u(t),v(t)}\colon{\Omega}_{u(t),v(t)}\to{\Omega}$ given by
$$T(x,z)=(x',z')=\left(x,\frac{z-v(t,x)}{u(t,x)-v(t,x)}\right).$$
With the definition of $\Omega_{u(t),v(t)}$ in Section~\ref{sec_intro}, it is easily checked that $T_{u(t),v(t)}$ is a diffeomorphism ${\Omega}_{u(t),v(t)}\to{\Omega}$ with the inverse
$$T^{-1}(x',z')=(x',z'(u(t,x')-v(t,x'))+v(t,x'))$$
and it is clear that $T$ and $T^{-1}$ can be extended to the boundary of ${\Omega}_{u(t),v(t)}$ and $\Omega$ respectively.
We introduce pull-back and push-forward operators $\theta^*(u,v)$ and $\theta_*(u,v)$ defined by $\theta^*(u,v)\tilde w=\tilde w\circ T_{u,v}$ and $\theta_*(u,v)w=w\circ T^{-1}_{u,v}$ where $w$ and $\tilde w$ are functions of the coordinates $(x,z)$ and $(x',z')$ respectively, i.e.,
$$[\theta^*(u,v)\tilde w](x,z)=\tilde w(T_{u,v}(x,z))\quad\text{and}\quad [\theta_*(u,v) w](x',z')=w(T_{u,v}^{-1}(x',z')).$$
Let $\widetilde{\Delta}_{u,v;\eps}=\theta_*(u,v)\Delta_\eps\theta^*(u,v)$ denote the time-dependent transformed Laplace operator on $\Omega$ which is explicitly given by
\begin{align}\widetilde{\Delta}_{u,v;\eps}\tilde w & = \eps^2\tilde w_{x'x'} - 2\eps^2\tilde w_{x'z'}\frac{z'(u_{x'}-v_{x'})+v_{x'}}{u-v}
+\tilde w_{z'z'}\frac{1+\eps^2[z'(u_{x'}-v_{x'})+v_{x'}]^2}{(u-v)^2} \nonumber\\
& \quad +\eps^2\tilde w_{z'}\left(2\frac{u_{x'}-v_{x'}}{(u-v)^2}[z'(u_{x'}-v_{x'})+v_{x'}]-\frac{z'(u_{x'x'}-v_{x'x'})+v_{x'x'}}{u-v}\right);\nonumber
\end{align}
here the notation $u_{x'}$ stands for $\partial_{x'}u$ et cetera.

For $q\in(2,\infty)$, we introduce the function spaces
$$
W^{2\alpha}_{q,D}(I):=\left\{
\begin{array}{ll}
\{w\in W^{2\alpha}_{q}(I);\,w(\pm 1)=0\}, & 2\alpha\in(1/q,2], \\
W_q^{2\alpha}(I), & 0\leq 2\alpha<1/q,
\end{array}
\right.
$$
and we define
$$
W^\alpha_{2,D}(\Omega):=\left\{
\begin{array}{ll}
\{w\in W^{\alpha}_2(\Omega);\,u|_{\partial\Omega}=0\}, & \alpha>1/2, \\
W^\alpha_{2}(\Omega), & 0\leq\alpha<1/2;
\end{array}
\right.
$$
the index $D$ indicates the Dirichlet boundary condition. The space $W^1_{2,D}(\Omega)$ is equipped with the norm $\norm{w}_{W^1_{2,D}(\Omega)}=\norm{\nabla w}_{L_2(\Omega)}$ and we will use the notation $W^{-\alpha}_{2,D}(\Omega)$ for the dual space $(W^{\alpha}_{2,D}(\Omega))'$, for $0\leq\alpha\leq 1$.

For $q\in(2,\infty)$ and $\kappa\in(0,1/2)$, we define the sets
\begin{align}
S_q(\kappa)&:=\bigg\{(u,v)\in W^2_{q}(I)\times W^2_{q}(I);\,(u,v)(\pm 1)=(0,-1),\,\norm{u}_{W^2_{q,D}(I)}<\frac{1}{\kappa},\nonumber\\
& \hspace{1cm}\norm{v+1}_{W^2_{q,D}(I)}<\frac{1}{\kappa},\, 2\kappa < u(x)-v(x) ,\,\forall x\in I\bigg\}\nonumber
\end{align}
and prepare the following lemma.
\lem
The sets $S_q(\kappa)+\{(0,1)\}\subset W^2_{q,D}(I)\times W^2_{q,D}(I)$ are open for $q\in(2,\infty)$ and $\kappa\in(0,1/2)$ and the closure of $S_q(\kappa)$ denoted as $\overline S_q(\kappa)$ is given by
\begin{align}
\overline S_q(\kappa)&=\bigg\{(u,v)\in W^2_{q}(I)\times W^2_{q}(I);\,(u,v)(\pm 1)=(0,-1),\,\norm{u}_{W^2_{q,D}(I)}\leq\frac{1}{\kappa},\nonumber\\
& \hspace{1cm}\norm{v+1}_{W^2_{q,D}(I)}\leq\frac{1}{\kappa},\, 2\kappa \leq u(x)-v(x) ,\,\forall x\in I\bigg\}.\label{closureS}
\end{align}
\endlem\rm
\proof For $q\in(2,\infty)$ and $\kappa\in(0,1/2)$ given, let $\tilde S_q(\kappa)$ be the set $S_q(\kappa)+\{(0,1)\}$. Then $(u,\tilde v)\in\tilde S_q(\kappa)$ if and only if $(u,\tilde v-1)\in S_q(\kappa)$ which is equivalent to
\begin{enumerate}
\item $u,\tilde v-1\in W^2_q(I)$,
\item $(u,\tilde v)(\pm 1)=(0,0)$,
\item $\norm{u}_{W^2_q(I)},\norm{\tilde v}_{W^2_q(I)}<1/\kappa$,
\item $2\kappa<u-\tilde v+1$ on $I$.
\end{enumerate}
The lemma claims that given $(u,\tilde v)\in \tilde S_q(\kappa)$ there exists $\eps>0$ such that
$$(u+\eps w_1,\tilde v+\eps w_2)\in\tilde S_q(\kappa)$$
for $w_1,w_2\in W^2_{q,D}(I)$, $\norm{w_1}_{W^2_q(I)},\norm{w_2}_{W^2_q(I)}<1$. To prove this, we first note that due to Sobolev's embedding theorem $W^2_q(I)\hookrightarrow C^1(I)$, so that there is a constant $c>0$ only depending on $q$ such that
$$\norm{w_1}_{\infty},\norm{w_{1x}}_{\infty},\norm{w_2}_{\infty},\norm{w_{2x}}_{\infty}\leq c.$$
Now we observe:
\begin{enumerate}
\item As $W^2_q(I)$ is a vector space, it is clear that $u+\eps w_1$ and $\tilde v+\eps w_2-1=(\tilde v-1)+\eps w_2$ belong to $W^2_q(I)$ for $(u,\tilde v-1)\in W^2_q(I)$ and $w_1,w_2\in W^2_{q,D}(I)$.
\item As $w_1,w_2$ have Dirichlet boundary conditions on $[-1,1]$ it is also clear that $(u+\eps w_1,\tilde v+\eps w_2)(\pm 1)=(0,0)$.
\item There exists $\delta_1>0$ such that $\norm{u}_{W^2_q(I)}\leq 1/\kappa-\delta_1$. Then
$$\norm{u+\eps w_1}_{W^2_q(I)}< 1/\kappa+\eps-\delta_1<1/\kappa$$
for $\eps<\delta_1$. The corresponding estimate for $\tilde v$ is obtained similarly for $\eps$ smaller than a number $\delta_2>0$.
\item There exists $\delta_3>0$ such that $u-\tilde v+1\geq 2\kappa+\delta_3$ on $I$. Then
\begin{align}
u+\eps w_1-(\tilde v+\eps w_2)+1 & = u-\tilde v+1+\eps(w_1-w_2)\nonumber\\
& \geq 2\kappa+\delta_3-2c\eps\nonumber\\
& >2\kappa\nonumber
\end{align}
for $\eps<\delta_3/(2c)$.
%
\end{enumerate}
Finally taking $\eps$ to be smaller than $\min\{\delta_1,\delta_2,\delta_3/(2c)\}$ achieves the proof of the first statement. Given a sequence $(u_n,v_n)\in S_q(\kappa)$ that converges to some $(u,v)\in W^2_q(I)\times W^2_q(I)$, it is immediately clear that $(u,v)$ belongs to the set on the right-hand side of \eqref{closureS}. It is an elementary proof that replacing some $<$ by $\leq$ in the definition of $S_q(\kappa)$ yields a subsets of the closure $\overline S_q(\kappa)$ with respect to $W^2_q(I)\times W^2_q(I)$. This achieves the proof of the lemma.
\endproof
\subsection{The elliptic problem}
We let $\tilde\phi(t,x',z')=\theta_*(u(t),v(t))\phi$, $\psi(t,x',z')=\tilde\phi(t,x',z')-z'$ and
$$f_{u,v;\eps} = \widetilde{\Delta}_{u,v;\eps}z' = \eps^2\left(2\frac{u_{x'}-v_{x'}}{(u-v)^2}[z'(u_{x'}-v_{x'})+v_{x'}]-\frac{z'(u_{x'x'}-v_{x'x'})+v_{x'x'}}{u-v}\right)$$
and rewrite the elliptic problem \eqref{originalproblem1}--\eqref{originalproblem2} as
\begin{align}
-\left(\widetilde{\Delta}_{u(t),v(t);\eps}\psi\right)(t,x',z') & = f_{u(t),v(t);\eps}, & \hspace{-1cm} (x',z')\in\Omega,\, t>0,\nonumber\\
\psi (t,x',z') & = 0, & \hspace{-1cm} (x',z')\in\partial\Omega,\, t>0.\nonumber
\end{align}
For $q\in(2,\infty)$, $\kappa\in(0,1/2)$ and $(u,v)\in\overline S_q(\kappa)$, the operator $-\widetilde{\Delta}_{u,v;\eps}$ is elliptic with an ellipticity constant independent of $(u,v)$ and we have the following lemma which generalizes \cite[Lemma~2.2]{ELW13}.
\lem\label{lem2.2} For each $(u,v)\in\overline S_q(\kappa)$ and $F\in W^{-1}_{2,D}(\Omega)$, there is a unique solution $\Phi\in W^1_{2,D}(\Omega)$ to the boundary value problem
\begin{align}
-\widetilde{\Delta}_{u,v;\eps}\Phi & = F, & \hspace{-2cm} \text{\rm in }\Omega,\nonumber\\
\Phi & = 0, & \hspace{-2cm} \text{\rm on }\partial\Omega,\nonumber
\end{align}
and there is a constant $C_1>0$ only depending on $\kappa$ and $\eps$ such that
$$\norm{\Phi}_{W^1_{2,D}(\Omega)}\leq C_1\norm{F}_{W^{-1}_{2,D}(\Omega)}.$$
Furthermore, if $F\in L_2(\Omega)$, then $\Phi\in W^2_{2,D}(\Omega)$ and
$$\norm{\Phi}_{W^2_{2,D}(\Omega)}\leq C_1\norm{F}_{L_2(\Omega)}.$$
\endlem\rm
\proof 
A careful observation shows that it suffices to establish the existence of positive constants $c_1,c_2$, only depending on $\kappa$ and $\eps$, such that
\beq\norm{\Phi}_{W^1_{2,D}(\Omega)}\leq c_1(\kappa,\eps)\norm{\Phi}_{L_2(\Omega)}+c_2(\kappa,\eps)\norm{F}_{W^{-1}_{2,D}(\Omega)}\label{eqc1c2}\eeq
for any test function $\Phi$ in the weak formulation of the Dirichlet problem $-\widetilde{\Delta}_{u,v;\eps}\Phi  = F$ on $\Omega$. Using the divergence form of $-\widetilde{\Delta}_{u,v;\eps}$, integration by parts and the Dirichlet boundary condition for $\Phi$, we obtain
\begin{align}
\ska{F}{\Phi} & = \int_\Omega\left[\eps^2\left(\Phi_{x'}-\frac{z'(u_{x'}-v_{x'})+v_{x'}}{u-v}\Phi_{z'}\right)^2+\frac{\Phi_{z'}^2}{(u-v)^2}\right]dx'dz' \nonumber\\
& \qquad + \eps^2\int_\Omega\frac{u_{x'}-v_{x'}}{u-v}\left[\frac{z'(u_{x'}-v_{x'})+v_{x'}}{u-v}\Phi_{z'}-\Phi_{x'}\right]\Phi\dx'dz'\nonumber
\end{align}
so that, setting
$$\xi=\frac{z'(u_{x'}-v_{x'})+v_{x'}}{u-v},$$
we get
\begin{align}
&\int_\Omega\left[\eps^2\left(\Phi_{x'}-\xi\Phi_{z'}\right)^2+\frac{\Phi_{z'}^2}{(u-v)^2}\right]dx'dz'
\leq\norm{F}_{W^{-1}_{2,D}(\Omega)}\norm{\Phi}_{W^1_{2,D}(\Omega)}\label{eqdualpair}\\
&\qquad+\eps^2\norm{\frac{u_{x'}-v_{x'}}{u-v}}_{L_\infty(I)}
\norm{\Phi_{x'}-\xi\Phi_{z'}}_{L_2(\Omega)}\norm{\Phi}_{L_2(\Omega)}.\nonumber
\end{align}
Using $u-v\geq 2\kappa$ and $\norm{u}_{C^1([-1,1])},\norm{v}_{C^1([-1,1])}\leq c_0(q,\kappa)$, an elementary computation shows that there exists a constant $0<\nu(\kappa,\eps)<1/2$ such that for any $(\zeta_1,\zeta_2)\in\R^2$
\beq\nu(\kappa,\eps)(\zeta_1^2+\zeta_2^2)\leq\eps^2\left(\zeta_1-\xi\zeta_2\right)^2
+\frac{\zeta_2^2}{(u-v)^2}.\label{nuabsch}\eeq
Letting $\zeta_1=\Phi_{x'}(x',z')$ and $\zeta_2=\Phi_{z'}(x',z')$ in \eqref{nuabsch} and integrating the inequality over $\Omega$, we can apply the resulting estimate twice to deduce from \eqref{eqdualpair} that \eqref{eqc1c2} holds true with $c_1=\eps c_0/(\kappa\sqrt\nu)$ and $c_2=\nu^{-1}$. The existence of a unique solution $\Phi\in W^1_{2,D}(\Omega)$ satisfying the estimates stated in the lemma now follows analogously to what has been done in \cite[Lemma~2.2]{ELW13}.
\endproof
An immediate consequence of Lemma~\ref{lem2.2} is that the transformed problem \eqref{originalproblem1}--\eqref{originalproblem2} on the fixed domain $\Omega$ has a unique solution $\tilde\phi_{u,v;\eps}\in W^2_2(\Omega)$ satisfying
\begin{align}
-\left(\widetilde{\Delta}_{u,v;\eps}\tilde\phi\right)(x',z') & = 0, & \hspace{-2cm} (x',z')\in\Omega,\label{phieq1trsf}\\
\tilde\phi (x',z') & = z', & \hspace{-2cm} (x',z')\in\partial\Omega.\label{phieq2trsf}
\end{align}
It is clear that, with the definition $(\tilde u,\tilde v)(x)=(u,v)(-x)$, $x\in I$, we have that $\tilde\phi_{\tilde u,\tilde v;\eps}(x',z')=\tilde\phi_{u,v;\eps}(-x',z')$, $(x',z')\in\Omega$.

Henceforth, we fix $\eps>0$ and omit it as an index to simplify notation. For $(u,v)\in\overline S_q(\kappa)$, let us define a second order linear operator $\mathcal A(u,v)\in\mathcal L(W^1_{2,D}(\Omega),W^{-1}_{2,D}(\Omega))$ by setting
$$\mathcal A(u,v)\Phi=-\widetilde\Delta_{u,v}\Phi,\quad\Phi\in W^1_{2,D}(\Omega).$$
A further consequence of Lemma~\ref{lem2.2} is that $\mathcal A(u,v)$ is invertible and it follows from the same arguments as in \cite[Lemma~2.3]{ELW13} that, for all $\theta\in[0,1]\backslash\{1/2\}$,
\beq\norm{\mathcal A(u,v)^{-1}}_{\mathcal L(W_{2,D}^{\theta-1},W^{\theta+1}_{2,D}(\Omega))}\leq C_2(\kappa,\eps),\quad\forall(u,v)\in\overline S_q(\kappa).\label{A-1thetaabsch}\eeq
We now show that $\tilde\phi_{u,v}$ depends Lipschitz continuously on $(u,v)\in \overline S_q(\kappa)$ in a suitable topology.
%
\lem\label{lemlipphi} Given $\xi\in[0,(q-1)/q)$ and $\alpha\in(\xi,1)$ there exists $C_3=C_3(\kappa,\eps)>0$ so that, for all $(u_1,v_1),(u_2,v_2)\in\overline S_q(\kappa)$,
$$\norm{\tilde\phi_{u_1,v_1}-\tilde\phi_{u_2,v_2}}_{W^{2-\alpha}_{2,D}(\Omega)}\leq C_3\norm{(u_1,v_1)-(u_2,v_2)}_{W^{2-\xi}_q(I)\times W^{2-\xi}_q(I)}.$$
\endlem\rm
\proof Since $$\tilde\phi_{u_1,v_1}-\tilde\phi_{u_2,v_2}=\psi_{u_1,v_1}-\psi_{u_2,v_2}=\mathcal A(u_1,v_1)^{-1}f_{u_1,v_1}-\mathcal A(u_2,v_2)^{-1}f_{u_2,v_2}$$ and $f_{u,v}\in L_2(\Omega)\hookrightarrow W^{-\alpha}_{2,D}(\Omega)$, for all $(u,v)\in\overline S_q(\kappa)$, the desired estimate follows immediately from the estimates
\begin{align}
\norm{\mathcal A(u_1,v_1)-\mathcal A(u_2,v_2)}_{\mathcal L(W^2_{2,D}(\Omega),W^{-\alpha}_{2,D}(\Omega))}&\leq c_1\norm{(u_1,v_1)-(u_2,v_2)}_{W^{2-\xi}_q(I)\times W^{2-\xi}_q(I)},\label{abschA}\\
\norm{f_{u_1,v_1}-f_{u_2,v_2}}_{W^{-\alpha}_{2,D}(\Omega)}&\leq c_2\norm{(u_1,v_1)-(u_2,v_2)}_{W^{2-\xi}_q(I)\times W^{2-\xi}_q(I)},\label{abschf}\\
\norm{\mathcal A(u_1,v_1)^{-1}}_{\mathcal L(W^{-\alpha}_{2,D}(\Omega),W^{2-\alpha}_{2,D}(\Omega))}&\leq c_3,\label{abschA-1}\\
\norm{\mathcal A(u_2,v_2)^{-1}}_{\mathcal L(L_2(\Omega),W^{2}_{2,D}(\Omega))}&\leq c_4,\label{abschA-1nr2}
\end{align}
where $c_1,\ldots,c_4$ are positive constants depending only on $\kappa$ and $\eps$. Note that \eqref{abschA-1nr2} is a direct consequence of Lemma~\ref{lem2.2} and \eqref{abschA-1} follows from \eqref{A-1thetaabsch} with $\theta=1-\alpha$. To prove \eqref{abschA} and \eqref{abschf}, we introduce the difference terms
\begin{align}
\gamma_0&:=\frac{1}{u_1-v_1}-\frac{1}{u_2-v_2}=\frac{(u_2-u_1)-(v_2-v_1)}{(u_1-v_1)(u_2-v_2)},\nonumber\\
\gamma_1&:=\frac{z'(u_1'-v_1')+v_1'}{u_1-v_1}-\frac{z'(u_2'-v_2')+v_2'}{u_2-v_2},\nonumber\\
\gamma_2&:=\frac{1+\eps^2(z'(u_1'-v_1')+v_1')^2}{(u_1-v_1)^2}-\frac{1+\eps^2(z'(u_2'-v_2')+v_2')^2}{(u_2-v_2)^2},\nonumber\\
\gamma_3&:=\frac{u_1'-v_1'}{(u_1-v_1)^2}(z'(u_1'-v_1')+v_1')-\frac{u_2'-v_2'}{(u_2-v_2)^2}(z'(u_2'-v_2')+v_2'),\nonumber\\
\gamma_4&:=\frac{z'(u_1''-v_1'')+v_1''}{u_1-v_1}-\frac{z'(u_2''-v_2'')+v_2''}{u_2-v_2},\nonumber
\end{align}
where $u_1'$ stands for $u_{1{x'}}$ et cetera. Consider $\Phi\in W^2_{2,D}(\Omega)$ and recall that $\mathcal A(u,v)\Phi\in L_2(\Omega)\hookrightarrow W^{-\alpha}_{2,D}(\Omega)$, so that for $\Psi\in W^{\alpha}_{2,D}(\Omega)$ we observe that
\begin{align} &\int_\Omega\left[\mathcal A(u_1,v_1)-\mathcal A(u_2,v_2)\right]\Phi\Psi\dx'dz'=-2\eps^2\int_\Omega\gamma_1\Phi_{x'z'}\Psi\dx'dz'\label{A1-A2}\\
&\qquad+\int_\Omega\gamma_2\Phi_{z'z'}\Psi\dx'dz'+2\eps^2\int_\Omega\gamma_3\Phi_{z'}\Psi\dx'dz'-\eps^2\int_\Omega\gamma_4\Phi_{z'}\Psi\dx'dz'.\nonumber
\end{align}
Rewriting $\gamma_1$ as
$$
\gamma_1=z'\left[\frac{u'_{1}-u'_{2}}{u_1-v_1}+u'_{2}\gamma_0\right]
+(1-z')\left[\frac{v'_{1}-v'_{2}}{u_1-v_1} +v'_{2}\gamma_0\right]
$$
and using that $W^{2-\xi}_q(I)\hookrightarrow W^1_\infty(I)$, $W^{\alpha}_{2,D}(\Omega)\hookrightarrow L_2(\Omega)$ and that $(u_1,v_1),(u_2,v_2)\in\overline S_q(\kappa)$, one concludes that
\begin{align}
&\left|\int_\Omega\gamma_1\Phi_{x'z'}\Psi\dx'dz'\right|\leq\norm{\gamma_1}_{L_\infty(\Omega)}\norm{\Phi_{x'z'}}_{L_2(\Omega)}\norm{\Psi}_{L_{2}(\Omega)}
\label{estgamma1}\\
&\qquad\leq c_5(\kappa)\norm{(u_1,v_1)-(u_2,v_2)}_{W^{2-\xi}_q(I)\times W^{2-\xi}_q(I)}\norm{\Phi}_{W^2_{2,D}(\Omega)}\norm{\Psi}_{W^{\alpha}_{2,D}(\Omega)}.\nonumber
\end{align}
Rewriting $\gamma_2$ and $\gamma_3$ as
$$
\gamma_2=\left(\frac{1}{u_1-v_1}+\frac{1}{u_2-v_2}\right)\gamma_0
+\eps^2\left(\frac{z'(u_1'-v_1')+v_1'}{u_1-v_1}+\frac{z'(u_2'-v_2')+v_2'}{u_2-v_2}\right)\gamma_1
$$
and
\begin{align}
\gamma_3 & = z'\left(\frac{u_1'-v_1'}{u_1-v_1}+\frac{u_2'-v_2'}{u_2-v_2}\right)\left(\frac{u_1'-u_2'}{u_1-v_1}+\frac{v_2'-v_1'}{u_1-v_1}+(u_2'-v_2')\gamma_0\right)\nonumber\\
&\quad+\left(\frac{v_2'}{u_2-v_2}+\frac{v_1'}{u_1-v_1}\right)\left(\frac{v_2'-v_1'}{u_2-v_2}-v_1'\gamma_0\right)
+v_1'\frac{u_1'-u_2'}{(u_1-v_1)^2}+u_2'\frac{v_1'-v_2'}{(u_2-v_2)^2}\nonumber\\
&\quad+u_2'v_1'\left(\frac{1}{u_1-v_1}+\frac{1}{u_2-v_2}\right)\gamma_0\nonumber
\end{align}
it is clear that
\begin{align}
&\left|\int_\Omega\gamma_2\Phi_{z'z'}\Psi\dx'dz'\right|,\left|\int_\Omega\gamma_3\Phi_{z'}\Psi\dx'dz'\right|\label{estgamma23}\\
&\qquad\leq c_6(\kappa)(1+\eps^2)\norm{(u_1,v_1)-(u_2,v_2)}_{W^{2-\xi}_q(I)\times W^{2-\xi}_q(I)}\norm{\Phi}_{W^2_{2,D}(\Omega)}\norm{\Psi}_{W^{\alpha}_{2,D}(\Omega)}.\nonumber
\end{align}
Writing $\gamma_4$ in the form
$$\gamma_4=z'\left[\frac{u''_{1}-u''_{2}}{u_1-v_1}+u''_{2}\gamma_0\right]
+(1-z')\left[\frac{v''_{1}-v''_{2}}{u_1-v_1} +v''_{2}\gamma_0\right]$$
and applying the generalized H\"older inequality, the fourth integral in \eqref{A1-A2} can be estimated by
\begin{align}
\left|\int_\Omega\gamma_4\Phi_{z'}\Psi\dx'dz'\right|&\leq\left|\int_\Omega\partial_{x'}^2(u_1-u_2)\frac{z'\Phi_{z'}\Psi}{u_1-v_1}\dx'dz'\right|\nonumber\\
&\hspace{-1cm}+\left|\int_\Omega\partial_{x'}^2(v_1-v_2)\frac{(1-z')\Phi_{z'}\Psi}{u_1-v_1}\dx'dz'\right|\nonumber\\
&\hspace{-1cm}+\norm{\gamma_0}_{L_\infty(I)}\left(\norm{u_2''}_{L_q(I)}+\norm{v_2''}_{L_q(I)}\right)
\norm{\Phi_{z'}}_{L_{2q/(q-2)}(\Omega)}\norm{\Psi}_{L_2(\Omega)}.\nonumber
\end{align}
For $\xi\in[0,(q-1)/q)$ one has $(W_{q'}^\xi(I))'=W^{-\xi}_q(I)$ so that
$$\left|\int_\Omega\partial_{x'}^2(u_1-u_2)\frac{z'\Phi_{z'}\Psi}{u_1-v_1}\dx'dz'\right|\leq\norm{u_1-u_2}_{W^{2-\xi}_q(I)}
\norm{\frac{1}{u_1-v_1}\int_0^1z'\Phi_{z'}\Psi\dz'}_{W^\xi_{q'}(I)}.$$
As explained in the proof of \cite[Lemma~2.4]{ELW13}, the second factor is bounded by $\norm{\Phi_{z'}}_{W^1_{2,D}(\Omega)}\norm{\Psi}_{W^{\alpha}_{2,D}(\Omega)}$, up to a positive constant only depending on $\kappa$. Clearly, the same arguments apply to the integral involving the factor $\partial_{x'}^2(v_1-v_2)$. Using that $W^1_2(\Omega)\hookrightarrow L_{2q/(q-2)}(\Omega)$ we infer
\begin{align}
&\left|\int_\Omega\gamma_4\Phi_{z'}\Psi\dx'dz'\right|\label{estgamma4}\\
&\qquad\leq c_7(\kappa)\norm{(u_1,v_1)-(u_2,v_2)}_{W^{2-\xi}_q(I)\times W^{2-\xi}_q(I)}\norm{\Phi}_{W^2_{2,D}(\Omega)}
\norm{\Psi}_{W^\alpha_{2,D}(\Omega)}.\nonumber
\end{align}
Now estimate \eqref{abschA} follows from \eqref{A1-A2}--\eqref{estgamma4}. Analogously, one deduces from
$$\int_\Omega(f_{u_1,v_1}-f_{u_2,v_2})\Psi\dx'dz'=2\eps^2\int_\Omega\gamma_3\Psi\dx'dz'-\eps^2\int_\Omega\gamma_4\Psi\dx'dz'$$
and the estimates
$$\left|\int_\Omega\gamma_3\Psi\dx'dz'\right|\leq c_8(\kappa)\norm{(u_1,v_1)-(u_2,v_2)}_{W^{2-\xi}_q(I)\times W^{2-\xi}_q(I)}\norm{\Psi}_{W^{\alpha}_{2,D}(\Omega)}$$
and
\begin{align}
\left|\int_\Omega\gamma_4\Psi\dx'dz'\right|&\leq\norm{u_1-u_2}_{W^{2-\xi}_q(I)}
\norm{\frac{1}{u_1-v_1}\int_0^1z'\Psi\dz'}_{W^\xi_{q'}(I)}\nonumber\\
&\hspace{-1cm}+\norm{v_1-v_2}_{W^{2-\xi}_q(I)}\norm{\frac{1}{u_1-v_1}\int_0^1(1-z')\Psi\dz'}_{W^\xi_{q'}(I)}\nonumber\\
&\hspace{-1cm}+\norm{\gamma_0}_{L_\infty(I)}\left(\norm{u_2''}_{L_2(I)}+\norm{v_2''}_{L_2(I)}\right)\norm{\Psi}_{L_2(\Omega)},\label{gamma4est}
\end{align}
applying once more the technique of \cite[Lemma~2.4]{ELW13} for the second factors of the first and second term on the right-hand side of \eqref{gamma4est}, that \eqref{abschf} holds true. This completes the proof of the lemma.
\endproof
A similar result with $\xi=0$ and $\alpha=0$ in the above lemma has been obtained in \cite{mk13_2}. In the following lemma, we show that the transformed right-hand sides of \eqref{originalproblem3}--\eqref{originalproblem4} depend analytically and Lipschitz continuously on $(u,v)\in\overline S_q(\kappa)$. To simplify notation, we write $u_x$ instead of $u_{x'}$ henceforth.
%
\lem\label{lemlipg} Let $q\in(2,\infty)$, $\kappa\in(0,1/2)$, $\eps>0$, $2\sigma\in[0,1/2)$ and $(u,v)\in\overline S_q(\kappa)$. Let $\tilde\phi_{u,v;\eps}\in W_2^2(\Omega)$ be the associated unique solution to \eqref{phieq1trsf}--\eqref{phieq2trsf}. Then the mapping $g_\eps\colon S_q(\kappa)\to W_{2,D}^{2\sigma}(I)\times W_{2,D}^{2\sigma}(I)$ defined by
$$
g_\eps(u,v)=\left(\frac{1+\eps^2u_{x}^2}{(u-v)^2}|\partial_{z'}\tilde\phi_{u,v;\eps}(\cdot,1)|^2,
\frac{1+\eps^2v_{x}^2}{(u-v)^2}|\partial_{z'}\tilde\phi_{u,v;\eps}(\cdot,0)|^2\right)
$$
is analytic, bounded, $g_\eps(0,-1)=(1,1)$, and if $\xi\in[0,1/2)$ and $\nu\in[0,(1-2\xi)/2)$, then there exists a constant $C_4(\kappa,\eps)>0$ such that
\beq\label{Lipg}\norm{g_\eps(u_1,v_1)-g_\eps(u_2,v_2)}_{W^\nu_2(I)\times W^\nu_2(I)}
\leq C_4\norm{(u_1,v_1)-(u_2,v_2)}_{W^{2-\xi}_q(I)\times W^{2-\xi}_q(I)}.\eeq
\endlem\rm
\proof We first recall from Lemma 8 of \cite{mk13_2} and the proof of Proposition 1 of \cite{mk13_2} that, for any $(u,v)\in\overline S_q(\kappa)$,
\begin{align}&\norm{\partial_{z'}\tilde\phi_{u,v}(\cdot,1)}_{W^{1/2}_2(I)}+\norm{\partial_{z'}\tilde\phi_{u,v}(\cdot,0)}_{W^{1/2}_2(I)}\nonumber\\
&\qquad+\norm{|\partial_{z'}\tilde\phi_{u,v}(\cdot,1)|^2}_{W^{2\sigma}_2(I)}+\norm{|\partial_{z'}\tilde\phi_{u,v}(\cdot,0)|^2}_{W^{2\sigma}_2(I)}\leq c_1(\kappa,\eps).\nonumber\end{align}
We rewrite $g_{\eps,1}(u_1,v_1)-g_{\eps,1}(u_2,v_2)$ as the sum of the three terms
\begin{align}
I_1&=\frac{1+\eps^2(u_1')^2}{(u_1-v_1)^2}\left(\partial_{z'}\tilde\phi_{u_1,v_1}(\cdot,1)+\partial_{z'}\tilde\phi_{u_2,v_2}(\cdot,1)\right)
\left(\partial_{z'}\tilde\phi_{u_1,v_1}(\cdot,1)-\partial_{z'}\tilde\phi_{u_2,v_2}(\cdot,1)\right),\nonumber\\
I_2&=(1+\eps^2(u_1')^2)|\partial_{z'}\tilde\phi_{u_2,v_2}(\cdot,1)|^2\frac{u_2-v_2+u_1-v_1}{(u_1-v_1)^2(u_2-v_2)^2}(u_2-u_1+v_1-v_2),\nonumber\\
I_3&=\eps^2|\partial_{z'}\tilde\phi_{u_2,v_2}(\cdot,1)|^2\frac{u_1'+u_2'}{(u_2-v_2)^2}(u_1'-u_2').\nonumber
\end{align}
For $2\sigma\in(\xi+\nu,1/2)$ and $s\in[\nu,1-\xi)$, $s\geq 1/q$, we have the continuous embeddings
$$W^s_q(I)\cdot W^{2\sigma}_2(I)\hookrightarrow W^\nu_2(I),\quad W^2_q(I)\cdot W^1_q(I)\cdot W^{1-\xi}_q(I)\hookrightarrow W^s_q(I)$$
and since $W^2_q(I)$ is an algebra, it is clear that
$$\norm{I_2}_{W^\nu_2(I)},\norm{I_3}_{W^\nu_2(I)}\leq c_2(\kappa,\eps)\norm{(u_1,v_1)-(u_2,v_2)}_{W^{2-\xi}_q(I)\times W^{2-\xi}_q(I)}.$$
By $W^1_q(I)\cdot W^{1/2}_2(I)\cdot W^{1/2-\alpha}_2(I)\hookrightarrow W^\nu_2(I)$, the algebra property of $W^1_q(I)$ and the regularity properties of the trace operator, cf.~\cite[Theorem 1.5.1.1]{G85}, we get
\begin{align}
\norm{I_1}_{W^\nu_2(I)}&\leq c_3(\kappa,\eps)\norm{\tilde\phi_{u_1,v_1}-\tilde\phi_{u_2,v_2}}_{W^{2-\alpha}_{2,D}(\Omega)}\nonumber\\
&\leq c_4(\kappa,\eps)\norm{(u_1,v_1)-(u_2,v_2)}_{W^{2-\xi}_q(I)\times W^{2-\xi}_q(I)},\nonumber
\end{align}
where we have used Lemma~\ref{lemlipphi}. The second component $g_{\eps,2}(u_1,v_1)-g_{\eps,2}(u_2,v_2)$ can be discussed similarly so that \eqref{Lipg} follows. Analyticity of the map $g_\eps$ follows from the analyticity of the maps $\mathcal A^{-1}\colon S_q(\kappa)\to\mathcal L(L_2(\Omega),W^2_{2,D}(\Omega))$ and $[(u,v)\mapsto f_{u,v}]$, $S_q(\kappa)\to L_2(\Omega)$. That $g_\eps(0,-1)=(1,1)$ and that $g_\eps$ is bounded is clear.
\endproof
\subsection{The abstract quasi-linear evolution equation}
Let $q\in(2,\infty)$, $\xi\in(0,\frac{q-1}{q})$ and $\kappa\in(0,1/2)$ and let $Z_q(\kappa)$ be the closed $1/\kappa$-ball in $W^{2-\xi}_q(I)$. We define, for $w_1\in Z_q(\kappa)$, the operator
\beq\label{defA}A(w_1)w_2:=-\frac{w_{2xx}}{(1+w_{1x}^2)^{3/2}},\quad D(A(w_1))=W^2_{q,D}(I).\eeq
Regarding \eqref{originalproblem6}, we introduce the function $\hat v=v+1$ and $\hat g_\eps(u,\hat v)=g_\eps(u,\hat v-1)=g_\eps(u,v)$ to rewrite the problem \eqref{originalproblem3}--\eqref{originalproblem8} as
\begin{align}
\frac{d}{dt}\begin{pmatrix} u\\\hat v \end{pmatrix} + \begin{pmatrix} A(\eps u) & 0 \\ 0 & A(\eps\hat v) \end{pmatrix}\begin{pmatrix} u\\\hat v \end{pmatrix}
&= \begin{pmatrix} -\lambda & 0 \\ 0 & \mu \end{pmatrix}\hat g_\eps(u,\hat v),&t>0,\label{uvprobtrsf1}\\
\begin{pmatrix} u\\\hat v \end{pmatrix}&=\begin{pmatrix} u_0\\\hat v_0 \end{pmatrix},&t=0.\label{uvprobtrsf2}
\end{align}
Note that the boundary conditions \eqref{originalproblem5}--\eqref{originalproblem6} are incorporated in the domain of the operator $A(\cdot)$. We now recall some important properties of $A(\cdot)$ from \cite{ELW13}:
For $\omega>0$ and $k\geq 1$ let $\mathcal H(W^2_{q,D}(I),L_q(I);k,\omega)$ be the set of all $A\in\mathcal L(W^2_{q,D}(I),L_q(I))$ such that $\omega+A$ is an isomorphism $W^2_{q,D}(I)\to L_q(I)$ satisfying
$$\frac{1}{k}\leq\frac{\norm{(\mu+A)z}_{L_q(I)}}{|\mu|\norm{z}_{L_q(I)}+\norm{z}_{W^2_{q,D}(I)}}\leq k,\quad\text{Re}(\mu)\geq\omega,\quad z\in W^2_{q,D}(I)\backslash\{0\}.$$
If $A\in\mathcal H(W^2_{q,D}(I),L_q(I);k,\omega)$, then $-A$ generates an analytic semigroup on $L_q(I)$ with domain $W^2_{q,D}(I)$. By \cite[Lemma~3.1]{ELW13}, for fixed $q\in (2,\infty)$, $\kappa\in(0,1/2)$ and $\xi\in(0,(q-1)/q)$, there are $k(\kappa)\geq 1$ and $\omega(\kappa)>0$ such that for any $w\in Z_q(\kappa)$, $-2\omega+A(w)\in \mathcal H(W^2_{q,D}(I),L_q(I);k,\omega)$ and $A(w)$ is resolvent positive satisfying
$$\norm{A(w_1)-A(w_2)}_{\mathcal L(W^2_{q,D}(I),L_q(I))}\leq\ell(\kappa)\norm{w_1-w_2}_{W^{2-\xi}_q(I)}$$
with a positive constant $\ell(\kappa)$. For $\rho\in(0,1)$ and $N,\tau>0$ let
\begin{align}
\mathcal W_\tau(\kappa)&:=\big\{w\in C([0,\tau],W^{2-\xi}_{q,D}(I));\,\norm{w(t)-w(s)}_{W^{2-\xi}_{q,D}(I)}\leq\frac{N}{\ell(\kappa)}|t-s|^\rho\nonumber\\
&\qquad \text{and }w(t)\in Z_q(\kappa)\text{ for }0\leq t,s\leq\tau\big\}.\nonumber
\end{align}
By \cite[Proposition~3.2]{ELW13}, there is a constant $c_*(\rho)>0$, independent of $N,\tau$, such that for each $w\in\mathcal W_\tau(\kappa)$ there exists a unique parabolic evolution operator $\mathcal U_{A(w)}(t,s)$, $0\leq s\leq t\leq\tau$, possessing $W^2_{q,D}(I)$ as a regular subspace and satisfying
$$\norm{\mathcal U_{A(w)}(t,s)}_{\mathcal L(W^{2\alpha}_{q,D}(I),W^{2\beta}_{q,D}(I))}\leq c_{**}(\kappa)(t-s)^{\alpha-\beta}e^{-\vartheta(t-s)},\quad 0\leq s<t\leq\tau,$$
for $0\leq\alpha\leq\beta\leq 1$ with $2\alpha,2\beta\neq 1/q$. The constant $c_{**}(\kappa)\geq 1$ depends on $N$, $\alpha$ and $\beta$ but is independent of $\tau$ and $-\vartheta=c_*(\rho)N^{1/\rho}-\omega(\kappa)$. Moreover $\mathcal U_{A(w)}(t,s)\in\mathcal L(L_q(I))$ is a positive operator for $0\leq s\leq t\leq\tau$.

Let $\lambda,\mu>0$, $q\in(2,\infty)$, $\eps\in(0,1)$ and fix $\kappa\in(0,1/4)$. As in the proof of \cite[Theorem~1.1]{ELW13}, we also fix $0<\xi<1/q$, $0<1/2-1/q<2\sigma<1/2-\xi$, $4\rho\in(0,\xi)$ and $N>0$ such that $-\vartheta<0$ and, for $w\in\mathcal W_\tau(\kappa)$ fixed,
\begin{align}
&\norm{\mathcal U_{A(w)}(t,s)}_{\mathcal L(W^2_{q,D}(I))}\label{abschU}\\
&\qquad+(t-s)^{-\sigma+1+\frac{1}{2}(\frac{1}{2}-\frac{1}{q})}\norm{\mathcal U_{A(w)}(t,s)}_{\mathcal L(W^{2\sigma}_{2,D}(I),W^{2}_{q,D}(I))}\leq M(\kappa)e^{-\vartheta(t-s)},\nonumber
\end{align}
for any $0\leq s\leq t\leq\tau$, with a constant $M(\kappa)\geq 1$ independent of $\omega$ and $\tau$. Note that \eqref{abschU} generalizes formula (34) of \cite{mk13_2}. We consider $(u_0,v_0)\in W^2_q(I)\times W^2_q(I)$ satisfying $(u_0,v_0)(\pm 1)=(0,-1)$, $-1\leq v_0<u_0\leq 0$ on $I$ and
$$\norm{u_0}_{W^2_{q,D}(I)},\norm{\hat v_0}_{W^2_{q,D}(I)}<\min\left\{\frac{1}{4}-\kappa,\frac{\kappa}{2M}\right\}=:r(\kappa).$$
In view of the continuous embedding $W^2_q(I)\hookrightarrow L_\infty(I)$ with embedding constant $2$, cf.\ the proof of Theorem~1.1 of \cite{ELW13}, this implies that $(u_0,v_0)\in S_q(2\kappa)$ and that $\norm{\mathcal U_{A(w)}(t,0)\hat v_0}_{L_\infty(I)}<\kappa$ on $[0,\tau]$. Increasing $M$ if necessary, we can also assume that $\norm{u_0}_{W^{2-\xi}_{q,D}(I)},\norm{\hat v_0}_{W^{2-\xi}_{q,D}(I)}\leq 1/\kappa$ without loss of generality.
%
%
Let $\kappa_0=\kappa/M\leq\kappa$ and define the spaces
\begin{align}
&X_\tau(\kappa):=\bigg\{(u,\hat v)\in\mathcal W_\tau(\kappa)\times\mathcal W_\tau(\kappa);\;\norm{u(t)}_{W^2_{q,D}(I)},\norm{\hat v(t)}_{W^2_{q,D}(I)}\leq\frac{1}{\kappa_0}\nonumber\\
&\hspace{2cm} 1+u(t)-\hat v(t)\geq 2\kappa,
\;\forall t\in[0,\tau]\bigg\}.\nonumber
\end{align}
Then $X_\tau(\kappa)\subset\overline S_q(\kappa_0)+\{(0,1)\}$, $\eps X_\tau(\kappa)\subset X_\tau(\kappa)$ and $X_\tau(\kappa)$ is a complete metric space with respect to the topology of $C([0,\tau];W^{2-\xi}_{q,D}(I)\times W^{2-\xi}_{q,D}(I))$. We now define, for $t\in[0,\tau]$ and $(u,\hat v)\in X_\tau(\kappa)$,
\begin{align}
F(u,\hat v)(t)&:=\begin{pmatrix}\mathcal U_{A(\eps u)}(t,0) & 0 \\ 0 & \mathcal U_{A(\eps\hat v)}(t,0)\end{pmatrix}\begin{pmatrix} u_0\\\hat v_0\end{pmatrix}\label{varofconstform}\\
&\qquad+\int_0^t\begin{pmatrix}-\lambda\,\mathcal U_{A(\eps u)}(t,s) & 0 \\ 0 & \mu\,\mathcal U_{A(\eps\hat v)}(t,s)\end{pmatrix}\hat g_\eps(u(s),\hat v(s))\,ds
\nonumber\end{align}
and claim that $F\colon X_\tau(\kappa)\to X_\tau(\kappa)$ is a contraction for either $(\lambda,\mu)$ arbitrary and $\tau$ sufficiently small or for $(\lambda,\mu)$ and $(u_0,\hat v_0)$ small and $\tau$ arbitrary. Recall from Lemma~\ref{lemlipg} that, for $(u_1,v_1),(u_2,v_2)\in\overline S_q(\kappa)$,
\begin{align}
&\norm{\hat g_\eps(u_1,\hat v_1)-\hat g_\eps(u_2,\hat v_2)}_{W^{2\sigma}_{2,D}(I)\times W^{2\sigma}_{2,D}(I)}\nonumber\\
&\qquad\leq C_4(\kappa,\eps)\norm{(u_1,v_1)-(u_2,v_2)}_{W^{2-\xi}_q(I)\times W^{2-\xi}_q(I)},\nonumber
\end{align}
and that
\beq\norm{\hat g_\eps(u,\hat v)}_{W^{2\sigma}_{2,D}(I)\times W^{2\sigma}_{2,D}(I)}\leq C_5(\kappa,\eps),\quad\forall(u,v)\in\overline S_q(\kappa).\label{boundg}\eeq
Let
$$\mathcal I(\tau):=\int_0^\tau e^{-\vartheta s}s^{\sigma-1-\tfrac{1}{2}(\tfrac{1}{2}-\tfrac{1}{q})}\, ds.$$
Then $\mathcal I\to 0$ as $\tau\to 0$, $\mathcal I\to\mathcal I(\infty)<\infty$ for $\tau\to\infty$ and $\tau\mapsto\mathcal I(\tau)$ is monotonically increasing on $[0,\infty)$. Using that $W^2_q(I)\hookrightarrow L_\infty(I)$ with embedding constant $2$ together with the positivity of the evolution operator and \eqref{boundg}, one concludes from \autoref{varofconstform} that, for $i=1,2$,
\begin{align}
0&\geq F_1(u,\hat v)(t),\label{boundsF1}\\
0&\leq F_2(u,\hat v)(t),\label{boundsF2}\\
1+F_1(u,\hat v)(t)-F_2(u,\hat v)(t)&\geq 4\kappa-\mathcal U_{A(\eps \hat v)}(t,0)\hat v_0-2(\lambda+\mu)M(\kappa)C_5(\kappa,\eps)\mathcal I(\tau)\text{ and}\label{diffF1F2}\\
\norm{F_i(u,\hat v)(t)}_{W^2_{q,D}(I)}&\leq\frac{1}{2\kappa_0}+\max\{\lambda,\mu\}M(\kappa)C_5(\kappa,\eps)\mathcal I(\tau).\label{normW2qD}
\end{align}
Applying \cite[II.\ Theorem~5.2.1]{A95} with $\alpha=1$, $\beta=1-\xi/2$ and $2\gamma=2\sigma-1/2+1/q$ together with
$$W^{2\sigma}_{2,D}(I)\hookrightarrow W_{q,D}^{2\sigma-\frac{1}{2}+\frac{1}{q}}(I)\hookrightarrow L_q(I)$$
we see that there exists a constant $C_6(\kappa)>0$ such that, with $m=\max\{\lambda,\mu\}$,
\begin{align}
&\norm{F(u_1,\hat v_1)(t)-F(u_2,\hat v_2)(t)}_{W^{2-\xi}_{q,D}(I)\times W^{2-\xi}_{q,D}(I)}\leq C_6\bigg(m\max_{0\leq t\leq\tau}\left(t^{\frac{\xi}{2}+\sigma-\tfrac{1}{2}(\tfrac{1}{2}-\tfrac{1}{q})}e^{-\vartheta t}\right)\label{LipF}\\
&\nonumber+\left(m+\frac{1}{2}\norm{(u_0,\hat v_0)}_{W^2_{q,D}(I)\times W^2_{q,D}(I)}\right)\max_{0\leq t\leq\tau}\left(t^{\frac{\xi}{2}}e^{-\vartheta t}\right)\bigg)
\norm{(u_1,\hat v_1)-(u_2,\hat v_2)}_{X_\tau(\kappa)}.
\end{align}
Applying \cite[II.\ Theorem~5.3.1]{A95} with $2\alpha=2-\xi+4\rho$ and $2\beta=2-\xi$ together with the embedding
$$W^{2}_{q,D}(I)\hookrightarrow W^{2-\xi+4\rho}_{q,D}(I)$$
we obtain, for $i=1,2$, $0\leq s\leq t\leq\tau$ and $(u,\hat v)\in X_\tau(\kappa)$,
\begin{align}
&\norm{F_i(u,\hat v)(t)-F_i(u,\hat v)(s)}_{W^{2-\xi}_{q,D}(I)}\leq C_7
\max_{0\leq t\leq\tau}\left(t^\rho e^{-\vartheta t}\right)\label{LipF2}\\
&\qquad\times\left(\norm{(u_0,\hat v_0)}_{W_{q,D}^{2-\xi+4\rho}(I)\times W_{q,D}^{2-\xi+4\rho}(I)}+2mC_5\right)(t-s)^\rho\nonumber
\end{align}
where $C_7(\kappa)>0$. As $F(u,\hat v)(0)=(u_0,\hat v_0)$, we conclude from \eqref{LipF2} and the triangle inequality that
\begin{align}
\norm{F_i(u,\hat v)(t)}_{W^{2-\xi}_{q,D}(I)}&\leq C_7
\max_{0\leq t\leq\tau}\left(t^{2\rho} e^{-\vartheta t}\right)\label{normboundF}\\
&\quad\times\left(\norm{(u_0,\hat v_0)}_{W_{q,D}^{2-\xi+4\rho}(I)\times W_{q,D}^{2-\xi+4\rho}(I)}+2mC_5\right)+\frac{1}{2\kappa}.\nonumber
\end{align}
It follows from \eqref{diffF1F2}--\eqref{normboundF} that we can choose $\tau>0$ sufficiently small so that $F\colon X_\tau(\kappa)\to X_\tau(\kappa)$ is indeed a contraction. The unique fixed point of $F$ in $X_\tau(\kappa)$ is a mild solution to \eqref{uvprobtrsf1}--\eqref{uvprobtrsf2} which can, according to \cite[Theorem~4.2]{A86} and \cite[Theorem 10.1]{A93}, be extended to a strong solution on a maximal interval of existence with the regularity specified in Theorem~\ref{thm_lwp}. Regarding \eqref{boundsF1}, \eqref{boundsF2}, Theorem~\ref{thm_lwp} follows from arguments very similar to what is presented in the proof of \cite[Theorem~2]{mk13_2} and \cite[Theorem~1.1]{ELW13}.
\subsection{The small aspect ratio limit}
We now establish that there is a positive $\eps$-independent lower bound for the maximal existence times $T_\eps$ of solutions $(u_\eps,v_\eps,\phi_\eps)$ to \eqref{originalproblem1}--\eqref{originalproblem8} as $\eps\to 0$. Then Theorem~\ref{thm_sar} follows from arguments very similar to what is presented in the proof of \cite[Theorem~10]{mk13_2} and \cite[Theorem~1.4]{ELW13}.

Fix $\lambda,\mu>0$, $q\in(2,\infty)$ and $\kappa\in(0,1/2)$ and consider $(u_0,v_0)\in W^2_{q}(I)\times W^2_{q}(I)$ with $(u_0,v_0)(\pm 1)=(0,-1)$, $-1\leq v_0<u_0\leq 0$ and $\norm{u_0}_{W^2_{q,D}(I)},\norm{v_0+1}_{W^2_{q,D}(I)}<r(\kappa/2)$
so that $(u_0,v_0)\in S_q(\kappa)$. For $\eps\in(0,1)$ we denote by $(u_\eps,v_\eps,\phi_\eps)$ the unique solution to \eqref{originalproblem1}--\eqref{originalproblem8} with initial values $(u_0,v_0)$, defined on the maximal interval $[0,T_\eps)$. Let $\kappa_1:=\kappa/(2M)<\kappa$ with $M$ as in \eqref{abschU} and
$$\tau_\eps:=\sup\left\{t\in[0,T_\eps);(u_\eps(s),v_\eps(s))\in\overline S_q(\kappa_1);\,\forall s\in[0,t]\right\}>0.$$
We then have $T_\eps\geq\tau_\eps$,
$$u_\eps(t)-v_\eps(t)\geq 2\kappa_1,\quad -1\leq v_\eps(t)<u_\eps(t)\leq 0\quad\text{on }[0,\tau_\eps]\times[-1,1],$$
and, by the continuous embedding $W_q^2(I)\hookrightarrow W^1_\infty(I)$,
$$\norm{u_\eps(t)}_{W_q^2(I)} + \norm{v_\eps(t)}_{W_q^2(I)} + \norm{u_\eps(t)}_{W_\infty^1(I)} + \norm{v_\eps(t)}_{W_\infty^1(I)}
\leq C_8(\kappa),\;\forall t\in[0,\tau_\eps].$$
Henceforth, we choose $\eps$ sufficiently small, precisely, $\eps$ smaller than some $\eps_*\in(0,1)$, so that
$$\eps_*^2\left(\norm{u_{\eps,x}(t)}_{L_\infty(I)}+2\norm{v_{\eps,x}(t)}_{L_\infty(I)}\right)^2\leq\frac{1}{2},
\quad\forall(t,\eps)\in[0,\tau_\eps]\times(0,\eps_*].$$
For $(t,x',z')\in[0,\tau_\eps]\times\overline\Omega$, we recall the definition $\psi_\eps(t,x',z')=\tilde\phi_\eps(t,x',z')-z'$, where
$\tilde\phi_\eps(t,x',z') = \theta_*(u(t),v(t))\phi_\eps$. Then $\psi_\eps(t)$ satisfies the uniform estimates established in \cite[Lemma~8]{mk13_2}. The fact that multiplication $W^1_q(I)\cdot W^{1/2}_2(I)\cdot W^{1/2}_2(I)\hookrightarrow W^{2\sigma}_2(I)$, $2\sigma\in(0,1/2)$, is continuous implies that
\beq\label{boundgohneeps}\norm{g_\eps(u_\eps(t),v_\eps(t))}_{W^{2\sigma}_2(I)\times W^{2\sigma}_2(I)}\leq C_9(\kappa).\eeq
Using \eqref{abschU}, \eqref{varofconstform}, \eqref{boundgohneeps} and that $(u_0,v_0)\in S_q(\kappa)$ we get
\beq\label{sarbounds}\norm{u_\eps(t)}_{W^2_{q,D}(I)},\norm{v_\eps(t)+1}_{W^2_{q,D}(I)}\leq\frac{M}{\kappa}+mMC_9\mathcal I(t),\eeq
with $m=\max\{\lambda,\mu\}$. Regarding \eqref{boundsF1}--\eqref{diffF1F2}, we recall that
\begin{align}
\label{boundsF1_}u_\eps(t)&\leq 0,\\
\label{boundsF2_}v_\eps(t)&\geq -1\text{ and}\\
\label{diffF1F2_}u_\eps(t)-v_\eps(t)&\geq 2\kappa-\mathcal U_{A(\eps (v+1))}(t,0)(v_0+1)-2(\lambda+\mu) MC_9\mathcal I(t).
\end{align}
As $\norm{v_0+1}_{W^2_{q,D}(I)}<\frac{2M-1}{4M^2}\kappa$, we have that $\norm{\mathcal U_{A(\eps (v+1))}(t,0)(v_0+1)}_{L_\infty(I)}<\kappa-\kappa_1$.
Furthermore, there exists $\tau>0$ such that $2(\lambda+\mu)MC_9\mathcal I(t)\leq\kappa-\kappa_1$ on $[0,\tau]$. Decreasing $\tau$ if necessary to guarantee that $mMC_9\mathcal I(t)\leq M/\kappa$, we conclude from \eqref{sarbounds} and \eqref{diffF1F2_} that $(u_\eps,v_\eps)(t)\in\overline S_q(\kappa_1)$ for all $t\in[0,\tau]\cap[0,\tau_\eps]$ and in particular $\tau_\eps\geq\tau$. With
$$\Lambda(\kappa):=\min\left\{\frac{1}{\kappa C_9\mathcal I(\infty)},\frac{(2M-1)\kappa}{8M^2C_9\mathcal I(\infty)}\right\}$$
it is also clear that, for $\lambda,\mu<\Lambda(\kappa)$, we obtain from \eqref{sarbounds} and \eqref{diffF1F2_} that $\tau_\eps\geq\tau$ for any $\tau>0$ and this implies that $T_\eps=\infty$.
\section{The non-existence of global solutions}\label{sec_gwp}
In this chapter, we first focus on the stationary version of \eqref{originalproblem1}--\eqref{originalproblem8}, i.e., the problem
\begin{align}
-\Delta_\eps\phi & = 0, & \text{in }\Omega_{u,v},\label{stationaryproblem1}\\
\phi & = \frac{z-v}{u-v}, & \text{on }\partial\Omega_{u,v},\label{stationaryproblem2}\\
u_{xx} & = \lambda(1+\eps^2u_x^2)^{5/2}|\phi_z(x,u(x))|^2,&x\in I,\label{stationaryproblem3}\\
v_{xx} & = -\mu(1+\eps^2v_x^2)^{5/2}|\phi_z(x,v(x))|^2,&x\in I,\label{stationaryproblem4}\\
u(\pm 1) & = 0,&\label{stationaryproblem5}\\
v(\pm 1) & =-1.&\label{stationaryproblem6}
\end{align}
Let us introduce the functions
$$J(r):=\int_0^r\frac{ds}{(1+s^2)^{5/2}}=\frac{r(2r^2+3)}{3(r^2+1)^{3/2}},\quad \tilde J(r):=J(-r)$$
and
$$\xi_0(\eps):=\min\left\{\frac{2J(\eps)}{\eps},\frac{2}{3\eps}\right\}.$$
Note that $J$ is strictly increasing, concave and maps $[0,\infty)$ to $[0,2/3)$.
It has been shown in the proof of \cite[Theorem~5]{mk13_2} that the potential satisfies
\beq\label{boundsphistat}z-v(x)\leq\phi(x,z)\leq 1+z-u(x),\quad\forall(x,z)\in\overline\Omega_{u,v}.\eeq
Using the upper bound in \eqref{boundsphistat} and the function $J$, the methods used in the proof of \cite[Theorem~1.3]{ELW13} imply that there can be no solution of \eqref{stationaryproblem3} with boundary condition \eqref{stationaryproblem5} provided $\lambda>\xi_0(\eps)$. We now make use of the lower bound for $\phi$ and infer from $\phi_z(x,v(x))\geq 1$ and \eqref{stationaryproblem4} that
$$\frac{v_{xx}}{(1+\eps^2v_x^2)^{5/2}}=-\frac{1}{\eps}\partial_x\tilde J(\eps v_x)\leq-\mu.$$
Without loss of generality, we assume that $v$ attains a maximum at $x_m\in(-1,0]$. Integrating the above inequality over $[x_m,x]$ for $x\in[0,1]$ implies that
$$\tilde J(\eps v_x)\geq\mu\eps x,\quad x\in[0,1].$$
Now either $\mu\eps\geq 2/3$ and then $\tilde J(\eps v_x(1))\geq 2/3$ which implies that $v_x(1)=-\infty$, a contradiction, or $\mu\eps<2/3$ and then, by Jensen's inequality,
$$\tilde J(-\eps v(0)-\eps)=\tilde J\left(\int_0^1\eps v_x\dx\right)\geq\int_0^1\tilde J(\eps v_x)\dx\geq\mu\frac{\eps}{2}.$$
If $\mu>2J(\eps)/\eps$, we obtain that $J(\eps v(0)+\eps)>J(\eps)$, i.e., $v(0)>0$, which is again a contradiction. This completes the proof of Theorem~\ref{thm_nonex1}.

We now present a proof of Theorem~\ref{thm_nonex2} and begin with the following lemma which refines the estimates \eqref{boundsphistat}. Recall that, by our assumptions, we concentrate on solutions $(u,v)$ to \eqref{originalproblem1}--\eqref{originalproblem8} such that $v,-u-1\leq c$, for some $c<0$, and that $(u,v)$ stays bounded in $W^2_q(I)\times W^2_q(I)$, i.e.,
%
%
by Sobolev's embedding theorem, $\norm{u}_{C^1([-1,1])}$ and $\norm{v}_{C^1([-1,1])}$ are bounded by a positive constant only depending on $q$.
\lem\label{lemmaxprinz} Let $\phi\in W^2_2(\Omega_{u,v})$ be a solution of $-\Delta_\eps\phi=0$ satisfying the boundary conditions $\phi(\pm 1,z)=1+z$, $\phi(x,u(x))=1$ and $\phi(x,v(x))=0$. Then there is $n\in 2\N$ such that, for all $(x,z)\in\overline\Omega_{u,v}$,
$$x^n+z\leq\phi(x,z)\leq 2+z-x^n.$$
\endlem\rm
\proof For some $n\in 2\N$, let $S^-_n(x,z)=x^n+z$ and $S^+_n(x,z)=2+z-x^n$. Then
$$-\Delta_\eps S_n^-=(-\eps^2\partial_x^2-\partial_z^2)S_n^-=-\eps^2n(n-1)x^{n-2}\leq 0$$
and we observe that $$S_n^-(\pm 1,z)= 1+z = \phi(\pm 1,z)$$ and $$S_n^-(x,u(x))= x^n+u(x) \leq 1 = \phi(x,u(x)).$$ As $v_x$ is uniformly bounded by a constant only depending on $q$ and by $v\leq c$, $c<0$, we shall make use of the fact that $x^n\to 0$, $n\to\infty$, pointwise in $I$, to obtain that $v(x)\leq -x^n$ or equivalently $$S^-_n(x,v(x))\leq\phi(x,v(x)),$$ for some $n\in 2\N$ and all $x\in I$.
%
%
As
$$-\Delta_\eps(S^-_n-\phi)\leq 0\,\text{ in }\Omega\quad\text{and}\quad(S^-_n-\phi)|_{\partial\Omega}\leq 0,$$
we can apply the weak maximum principle to conclude that $S_n^-\leq\phi$ in $\overline\Omega_{u,v}$. Similarly, one shows that
$$-\Delta_\eps(S^+_n-\phi)\geq 0\,\text{ in }\Omega\quad\text{and}\quad(S^+_n-\phi)|_{\partial\Omega}\geq 0$$
so that the weak maximum principle implies that $S_n^+\geq\phi$ in $\overline\Omega_{u,v}$.
\endproof
Note that the number $n$ in the above lemma only depends on $c$ and $q$. Let us now modify the calculations in \cite{ELW13b} for the problem under discussion.

We multiply $\eps^2\phi_{xx}+\phi_{zz}=0$ by the function $\phi_z-1$, integrate over $\Omega_{u,v}$ and use integration by parts to obtain that
\begin{align}
0 & = -\eps^2\int_{\Omega_{u,v}}\phi_{x}\phi_{xz}\dx\d z+\eps^2\int_{\partial\Omega_{u,v}}\phi_x(\phi_z-1)n_1\d s\nonumber\\
&\qquad + \int_{\Omega_{u,v}}\left(\frac{1}{2}\frac{d}{dz}\phi_z^2-\phi_{zz}\right)dx\d z,\nonumber
\end{align}
with $n=(n_1,n_2)$ denoting the outward normal of $\partial\Omega_{u,v}$. Using the identities
\begin{align}
\phi_x(x,u(x))&=-u_x\phi_z(x,u(x)),\label{chainruleu}\\
\phi_x(x,v(x))&=-v_x\phi_z(x,v(x)),\label{chainrulev}
\end{align}
which follow from differentiating the boundary conditions $\phi(x,u(x))=1$ and $\phi(x,v(x))=0$, and $\phi_z(\pm 1,z)=1$ we obtain
\begin{align}
0 & = -\frac{\eps^2}{2}\int_{\Omega_{u,v}}\frac{d}{dz}\phi_x^2\dx\d z+\eps^2\int_I\phi_z(x,u(x))(\phi_z(x,u(x))-1)u_x^2\dx\nonumber\\
& \quad -\eps^2\int_I\phi_z(x,v(x))(\phi_z(x,v(x))-1)v_x^2\dx + \frac{1}{2}\int_I\left(\phi_z^2(x,u(x))-\phi_z^2(x,v(x))\right)dx\nonumber\\
& \quad - \int_I(\phi_z(x,u(x))-\phi_z(x,v(x)))\dx \nonumber\\
& = \frac{\eps^2}{2}\int_{I}\left(\phi_z^2(x,u(x))u_x^2-\phi_z^2(x,v(x))v_x^2\right)dx + \frac{1}{2}\int_I\left(\phi_z^2(x,u(x))-\phi_z^2(x,v(x))\right)dx\nonumber\\
& \quad -\eps^2\int_I\left(\phi_z(x,u(x))u_x^2-\phi_z(x,v(x))v_x^2\right)dx - \int_I(\phi_z(x,u(x))-\phi_z(x,v(x)))\dx \nonumber
\end{align}
and thus
\begin{align}
\int_I(1+\eps^2v_x^2)\phi_z^2(x,v(x))\dx & = \int_I(1+\eps^2u_x^2)\left(\phi_z^2-2\phi_z\right)(x,u(x))\dx\nonumber\\
&\quad +2\int_I(1+\eps^2v_x^2)\phi_z(x,v(x))\dx\nonumber\\
&\geq 2\int_I(1+\eps^2v_x^2)\phi_z(x,v(x))\dx - \int_I(1+\eps^2u_x^2)\dx. \nonumber
\end{align}
As $\norm{u}_{C^1([-1,1])}$ is bounded by a constant only depending on $q$, there is $\eps_0>0$ such that, for all $\eps<\eps_0$, we have that $\eps^2\norm{u_x}_{L_\infty(I)}^2\leq\frac{1}{4}$ and thus
\beq\label{step1}\int_I(1+\eps^2v_x^2)\phi_z^2(x,v(x))\dx \geq 2\int_I(1+\eps^2v_x^2)\phi_z(x,v(x))\dx - \frac{5}{2}.\eeq
A corresponding estimate with $v$ replaced by $u$ can be obtained similarly. Working with \eqref{step1} henceforth motivates to assume that $\mu\geq\lambda$ in the following, without loss of generality.

We multiply $\eps^2\phi_{xx}+\phi_{zz}=0$ by the function $\phi-1$, integrate over $\Omega_{u,v}$ and use integration by parts, \eqref{chainrulev} and Theorem~\ref{thm_lwp}.(iii) to obtain that
\begin{align}
\int_{\Omega_{u,v}}(\eps^2\phi_x^2+\phi_z^2)\dx\d z &= \eps^2\int_{\partial\Omega_{u,v}}\phi_x(\phi-1)n_1\d s+\int_{\partial\Omega_{u,v}}\phi_z(\phi-1)n_2\d s
\nonumber\\
& = -\eps^2\int_I\phi_x(x,v(x))v_x\dx+\eps^2\int_{-1}^0\phi_x(1,z)z\d z\nonumber\\
& \quad-\eps^2\int_{-1}^0\phi_x(-1,z)z\d z+\int_I\phi_z(x,v(x))\dx\nonumber\\
& = \int_I(1+\eps^2v_x^2)\phi_z(x,v(x))\dx + 2\eps^2\int_{-1}^0\phi_x(1,z)z\dz.\nonumber
\end{align}
By Lemma~\ref{lemmaxprinz},
$$\phi(x,z)-\phi(1,z)\leq 2+z-x^n-(1+z)=-(x^n-1)$$
and, for $x<1$, we obtain
$$\frac{\phi(x,z)-\phi(1,z)}{x-1}\geq -x^{n-1}-x^{n-2}-\ldots-x-1.$$
Sending $x\to 1$, we conclude that $\phi_x(1,z)\geq -n$ and thus, as $z\in[-1,0]$,
$$\int_{-1}^0\phi_x(1,z)z\d z\leq\frac{n}{2}.$$
This yields
\beq\label{step2}\int_{\Omega_{u,v}}(\eps^2\phi_x^2+\phi_z^2)\dx\d z\leq \int_I(1+\eps^2v_x^2)\phi_z(x,v(x))\dx +\eps^2n.\eeq
Now
\begin{align}\frac{1}{u-v}&=\frac{(\phi(x,u(x))-\phi(x,v(x)))^2}{u-v}=\frac{1}{u-v}\left(\int_{v(x)}^{u(x)}\phi_z(x,z)\d z\right)^2\label{step3}\\
&\leq\int_{\Omega_{u,v}}(\eps^2\phi_x^2+\phi_z^2)\dx\d z.\nonumber\end{align}
The function $\alpha(r):=\frac{1}{1+r}$, $r\in(-1,\infty)$ is convex and Jensen's inequality implies
$$\frac{1}{2}\int_I\frac{1}{u-v}\dx=\frac{1}{2}\int_I\frac{1}{1+[u-(v+1)]}\dx\geq\frac{1}{1+\frac{1}{2}\int_I[u-(v+1)]\dx}.$$
Setting
$$E(t):=-\frac{1}{2}\int_I[u-(v+1)]\dx,$$
we note that $E(t)\in[0,1)$, and using \eqref{step1}, \eqref{step2} and \eqref{step3}, we derive the inequality
\begin{align}
\frac{1}{1-E(t)}&\leq\frac{1}{2}\int_{\Omega_{u,v}}(\eps^2\phi_x^2+\phi_z^2)\dx\d z\nonumber\\
&\leq\frac{1}{2}\int_I(1+\eps^2v_x^2)\phi_z(x,v(x))\dx +\eps^2\frac{n}{2}\nonumber\\
&\leq\frac{1}{4}\left(\int_I(1+\eps^2v_x^2)\phi_z^2(x,v(x))\dx +\frac{5}{2} \right)+\eps^2\frac{n}{2}.\nonumber
\end{align}
In view of \eqref{originalproblem3} and \eqref{originalproblem4}, we observe that
\begin{align}
\frac{dE(t)}{dt}&=-\left[\frac{u_x}{2\sqrt{1+\eps^2u_x^2}}-\frac{v_x}{2\sqrt{1+\eps^2v_x^2}}\right]^{1}_{-1}
+\frac{\lambda}{2}\int_I\phi_z^2(x,u(x))(1+\eps^2u_x^2)\dx\nonumber\\
&\quad+\frac{\mu}{2}\int_I\phi_z^2(x,v(x))(1+\eps^2v_x^2)\dx\nonumber\\
&\geq-\frac{2}{\eps}+\frac{\mu}{2}\left(\frac{4}{1-E(t)}-\frac{5}{2}-2\eps^2n\right).\nonumber
\end{align}
If necessary, we decrease $\eps_0>0$ to guarantee that $2\eps^2n\leq\frac{1}{2}$ and so
$$\frac{dE(t)}{dt}\geq-\frac{2}{\eps}+2\mu\alpha(-E)-\frac{3}{2}\mu=:F_\mu(E).$$
As $-1<-E\leq 0$, $1\leq\alpha(-E)<\infty$ and hence $F_\mu(E)\geq -2/\eps+\mu/2$. If $\mu>4/\eps$, then $F_\mu(E)>0$ and the above inequality implies that $E(t)$ is strictly increasing. As $F_\mu(E)$ is also strictly increasing, we must have
$$\frac{dE(t)}{dt}\geq F_\mu(E(0))\geq F_\mu(0).$$
This shows that
$$1>E(t)\geq E(0)+F_\mu(0)t,\quad\forall t\in[0,T_\eps),$$
which immediately yields
$$T_\eps\leq T_\eps^*:=\frac{1-E(0)}{F_\mu(0)}=\frac{1}{\mu-4/\eps}\int_I(u_0-v_0)\dx$$
and hence \eqref{Tepsabsch}. Moreover, $0\leq\min_{x\in I}\{u(t)-v(t)\}\leq 1-E(t)$ so that, for $T_\eps=T_\eps^*$, $\liminf_{t\to T_\eps}\min_{x\in I}\{u(t)-v(t)\}=0$.
This completes the proof of Theorem~\ref{thm_nonex2}.
\section{Asymptotically stable steady state solutions}\label{sec_steady}
In terms of the coordinates $(x',z')\in\overline\Omega$, the problem \eqref{stationaryproblem1}--\eqref{stationaryproblem6} reads
\begin{align}
-\tilde\Delta_\eps\tilde\phi & = 0, & \text{in }\Omega,\label{stationaryproblem1'}\\
\tilde\phi & = z', & \text{on }\partial\Omega,\label{stationaryproblem2'}\\
u_{x'x'} & = \lambda\frac{(1+\eps^2u_{x'}^2)^{5/2}}{(u-v)^2}|\tilde\phi_{z'}(x',1)|^2,&x'\in I,\label{stationaryproblem3'}\\
v_{x'x'} & = -\mu\frac{(1+\eps^2v_{x'}^2)^{5/2}}{(u-v)^2}|\tilde\phi_{z'}(x',0)|^2,&x'\in I,\label{stationaryproblem4'}\\
u(\pm 1) & = 0,&\label{stationaryproblem5'}\\
v(\pm 1) & =-1.&\label{stationaryproblem6'}
\end{align}
Fix $q\in(2,\infty)$ and $\kappa\in(0,1/2)$. We recall the notation $\hat v=v+1$ and the definition of the operator $A$ in \eqref{defA}. Defining $h_\eps:=(h_{1,\eps},h_{2,\eps})\colon S_q(\kappa)\to L_q(I)\times L_q(I)$ by
\begin{align}
h_{1,\eps}(u,v):=\frac{(1+\eps^2u_{x'}^2)^{5/2}}{(u-v)^2}|\tilde\phi_{z'}(x',1)|^2\nonumber\\
h_{2,\eps}(u,v):=\frac{(1+\eps^2v_{x'}^2)^{5/2}}{(u-v)^2}|\tilde\phi_{z'}(x',0)|^2\nonumber
\end{align}
and recalling that $-A(0)=\partial_{x'}^2\in\mathcal L(W^2_{q,D}(I),L_q(I))$ is invertible, we introduce a map $F\colon\R^2\times S_q(\kappa)\to W^2_{q,D}(I)\times W^2_{q}(I)$ by
$$F(\Lambda,U):=\begin{pmatrix}U_1\\U_2\end{pmatrix}+\begin{pmatrix}\Lambda_1 & 0 \\ 0 & -\Lambda_2\end{pmatrix}A(0)^{-1}h_\eps(U_1,U_2).$$
Then $F(0,0)=(0,0)$ and $D_UF(0,0)=\text{id}$ so that, in view of the Implicit Function Theorem, there is $\delta=\delta(\kappa)>0$ and an analytic map $[\Lambda\mapsto U_\Lambda]\colon[0,\delta)^2\to W^2_{q,D}(I)\times W^2_{q}(I)$ satisfying $F(\Lambda,U_\Lambda)=0$. For $\Lambda\neq (0,0)$, let $\Phi_\Lambda$ be the potential associated with $U_\Lambda$. Then $(U_\Lambda,\Phi_\Lambda)\in S_q(\kappa)\times W^2_2(\Omega)$ is the unique stationary solution to \eqref{originalproblem1}--\eqref{originalproblem8}. Given $U=(U_1,U_2)$, we use the notation $\hat U=(U_1,\hat U_2)$ and we write $\Lambda=(\lambda,\mu)$. Letting $\hat V=U-U_\Lambda=\hat U-\hat U_\Lambda$ and introducing a map $Q=(Q_1,Q_2)$ by setting
\begin{align}
Q_1(u,\hat v)&:=-A(\eps u)u-\lambda\hat g_{\eps,1}(u,\hat v), \nonumber\\
Q_2(u,\hat v)&:=-A(\eps\hat v)\hat v+\mu\hat g_{\eps,2}(u,\hat v), \nonumber
\end{align}
we observe that $Q(\hat U_\Lambda)=0$, and we introduce the function
$$G_\Lambda(\hat V):=Q(\hat V+\hat U_\Lambda)-DQ(\hat U_\Lambda)\hat V$$
so that, for $\hat U$ being a solution of \eqref{originalproblem3}--\eqref{originalproblem4},
$$\frac{d}{dt}\hat V-DQ(\hat U_\Lambda)\hat V=G_\Lambda(\hat V).$$
Clearly, $G_\Lambda\in C^{\infty}(\mathcal O_\Lambda,L_q(I)\times L_q(I))$, where $\mathcal O_\Lambda\subset W^2_{q,D}(I)\times W^2_{q,D}(I)$ is a neighborhood of zero such that $U_\Lambda+\mathcal O_\Lambda\subset S_q(\kappa)$, $G_\Lambda(0)=0$ and $DG_\Lambda(0)=0$. A straightforward computation shows that
\begin{align}-DQ(\hat U_\Lambda)\hat V&=\begin{pmatrix}A(\eps U_{\Lambda,1}) & 0 \\ 0 & A(\eps\hat U_{\Lambda,2}) \end{pmatrix}\hat V \nonumber\\
& \quad + 3\eps^2\begin{pmatrix}\lambda g_{\eps,1}(U_\Lambda)\frac{\partial_{x'}U_{\Lambda,1}}{1+\eps^2(\partial_{x'}U_{\Lambda,1})^2} & 0 \\
0 & -\mu g_{\eps,2}(U_\Lambda)\frac{\partial_{x'}U_{\Lambda,2}}{1+\eps^2(\partial_{x'}U_{\Lambda,2})^2}\end{pmatrix}\partial_{x'}\hat V\nonumber\\
& \quad + \begin{pmatrix}\lambda & 0 \\ 0 & -\mu\end{pmatrix}D\hat g_\eps(\hat U_\Lambda)\hat V\nonumber\\
&=:\begin{pmatrix}A(\eps U_{\Lambda,1}) & 0 \\ 0 & A(\eps\hat U_{\Lambda,2}) \end{pmatrix}\hat V+B_\Lambda\hat V\nonumber
\end{align}
and we obtain that
$$\frac{d}{dt}\hat V+\left[\begin{pmatrix}A(\eps U_{\Lambda,1}) & 0 \\ 0 & A(\eps\hat U_{\Lambda,2}) \end{pmatrix}+B_\Lambda\right]
\hat V=G_\Lambda(\hat V).$$
Since $U_\Lambda\in S_q(\kappa)$, we have that
$$A(\eps U_{\Lambda,1}),A(\eps\hat U_{\Lambda,2})\in\mathcal H(W^2_{q,D}(I),L_q(I);k,\omega)$$ with a spectral bound less than $-\omega<0$. Since
$$\norm{B_\Lambda}_{\mathcal L(W^2_{q,D}(I)\times W^2_{q,D}(I),L_q(I)\times L_q(I))}\to 0,\quad\Lambda\to 0,$$
the operator $-(\text{diag}(A(\eps U_{\Lambda,1}),A(\eps\hat U_{\Lambda,2}))+B_\Lambda)$ generates an analytic semigroup on $L_q(I)\times L_q(I)$ with a negative spectral bound, cf.~\cite{ELW13} for more details in a similar case. Then we can apply \cite[Theorem~9.1.2]{Lunardi} to conclude Theorem~\ref{thm_stability}. From Theorem~\ref{thm_stability}.(ii) and the Lipschitz continuity of $\tilde\phi$ obtained in \cite[Proposition~1]{mk13_2}, we also conclude that
$$\norm{\tilde\phi_{u,v}-\Phi_\Lambda}_{W^2_2(\Omega)}\leq R'e^{-\omega_0 t}\norm{(u_0,v_0)-U_{\Lambda}}_{W_{q,D}^2(I)\times W_{q,D}^2(I)},\quad\forall t\geq 0.$$

\end{document}